\documentclass[a4paper,11pt]{amsart}
\usepackage{eurosym}
\usepackage{amsmath}
\usepackage{mathrsfs}
\usepackage{amsfonts}
\usepackage{graphicx}
\usepackage{color}
\usepackage{amsfonts}
\usepackage{amssymb}
\usepackage{amsthm}
\usepackage[hidelinks]{hyperref}
\usepackage[normalem]{ulem}
\usepackage{esint}
\usepackage{palatino}
\usepackage[abbrev,backrefs]{amsrefs}
\usepackage[hidelinks]{hyperref}
\usepackage{mathtools}
\usepackage{tikz-cd}
\usepackage{float}
\usepackage[text={6.8in,10in},centering]{geometry}
\newtheorem{theorem}{Theorem}[section]

\newtheorem{corollary}[theorem]{Corollary}

\newtheorem{lemma}[theorem]{Lemma}

\newtheorem{proposition}[theorem]{Proposition}
\newtheorem{remark}[theorem]{Remark}

\numberwithin{equation}{section}

\def\O{{\mathcal O}}

\def\dist{{\rm dist\,}}

\def\supp{{\rm supp\,}}

\newcommand{\enf}{\mathop{\mbox{$\essinf_\beta$}}}
\newcommand{\esnf}{\mathop{\mbox{$\essinf_{\beta_1}$}}}
\newcommand{\R}{\mathbb{R}}
\newcommand{\Rn}{\mathbb{R}^n}
\newcommand{\dx}{\,{\rm d}x}
\newcommand{\dy}{\,{\rm d}y}

\newcommand{\dt}{\,{\rm d}t}
\newcommand{\avgdb}[2]{f_{#1,#2}}
\newcommand{\avgd}[1]{f_{#1,\beta}}
\newcommand{\avgl}[1]{f_{#1}}
\newcommand{\cont}{\mathcal{H}_{\infty}^{\beta}}
\newcommand{\dc}{\,{\rm d}\cont}
\newcommand{\contq}{\mathcal{H}_{\infty}^{\beta,Q'}}
\newcommand{\dcq}{\,{\rm d}\contq}
\newcommand{\contb}[1]{\mathcal{H}_{\infty}^{#1}}
\newcommand{\dcb}[1]{\,{\rm d}\mathcal{H}_{\infty}^{#1}}
\newcommand{\cbmo}{\text{BMO}}
\newcommand{\bmo}{\text{BMO}^{\beta}}
\newcommand{\bmob}[1]{\text{BMO}^{#1}}
\newcommand{\bmobp}[2]{\text{BMO}^{#1,#2}}
\newcommand{\cblo}{\text{BLO}}
\newcommand{\blo}{\text{BLO}^{\beta}}
\newcommand{\bloq}{\text{BLO}^{\beta}_{*}}
\newcommand{\blob}[1]{\text{BLO}^{#1}}
\newcommand{\cvmo}{\text{VMO}}
\newcommand{\vmo}{\text{VMO}^{\beta}}
\newcommand{\vmob}{\text{VMO}^{\beta_1}}
\newcommand{\fma}{\mathcal{M}_{\alpha}}
\newcommand{\fdma}{\mathcal{M}_{\alpha}^{\beta}}

\newcommand{\fmal}{\mathcal{M}_{\alpha,\kappa}^{\text{loc}}}
\newcommand{\fmag}{\mathcal{M}_{\alpha,\kappa}^{\text{glob}}}
\newcommand{\maxl}{\mathcal{M}_{\alpha, \lambda r}^{\mathrm{loc}}}
\newcommand{\maxgg}{\mathcal{M}_{\alpha, \lambda r}^{\mathrm{glob}}}

\usepackage{scalerel,stackengine}
\stackMath

\DeclareMathOperator*{\essinf}{ess\,inf}
\DeclareMathOperator*{\argmin}{arg\,min}

\newcommand{\barint}{
	\rule[.036in]{.12in}{.009in}\kern-.16in \displaystyle\int }

\newcommand{\barcal}{\mbox{$ \rule[.036in]{.11in}{.007in}\kern-.128in\int $}}

\usepackage{enumitem}
\makeatletter
\let\@wraptoccontribs\wraptoccontribs
\makeatother

\mathchardef\mhyphen="2D

\title[Uncentered Fractional Maximal functions and mean oscillation spaces]{Uncentered Fractional Maximal functions and mean oscillation spaces associated with dyadic Hausdorff content}

\author[R. Basak]{Riju Basak}
	\address[R. Basak]{Department of Mathematics, National Taiwan Normal University, No. 88, Section 4, Tingzhou Road, Wenshan District, Taipei City, Taiwan 116, R.O.C.
	}
	\email{rijubasak52@ntnu.edu.tw}
	
	\author[Y.-W. B. Chen]{You-Wei Benson Chen}
	\address[Y.-W. B. Chen]{Department of Mathematics, National Taiwan Normal University, No. 88, Section 4, Tingzhou Road, Wenshan District, Taipei City, Taiwan 116, R.O.C.
	}
	\email{bensonchen.sc07@nycu.edu.tw}
	
	\author[P. Roychowdhury]{Prasun Roychowdhury}
	\address[P. Roychowdhury]{Dipartimento di Matematica e Applicazioni, Università degli Studi di Milano–Bicocca, Via Cozzi 55, 20125 Milano, Italy.
	} 
	\email{prasun.roychowdhury@unimib.it}
    
\date{\today}

\begin{document}

\begin{abstract}
We study the action of uncentered fractional maximal functions on mean oscillation spaces associated with the dyadic Hausdorff content $\mathcal{H}_{\infty}^{\beta}$ with $0<\beta\leq n$. For $0 < \alpha < n$, we refine existing results concerning the action of the Euclidean uncentered fractional maximal function $\mathcal{M}_{\alpha}$ on the functions of bounded mean oscillations (BMO) and vanishing mean oscillations (VMO). In addition, for $0 < \beta_1 \leq \beta_2 \leq n$, we establish the boundedness of the $\beta_2$-dimensional uncentered maximal function $\mathcal{M}^{\beta_2}$ on the space $\text{BMO}^{\beta_1}(\mathbb{R}^n)$, where $\text{BMO}^{\beta_1}(\mathbb{R}^n)$ denotes the mean oscillation space adapted to the dyadic Hausdorff content $\mathcal{H}_{\infty}^{\beta_1}$ on $\mathbb{R}^n$.
\end{abstract}

\maketitle
		
\section{Introduction and main results}
For $\alpha \in [0,n)$, the uncentered fractional maximal function $\mathcal{M}_\alpha$ is defined by 
\begin{align*} 
\mathcal{M}_\alpha f(x) := \sup_{Q \ni x} \frac{1}{\ell(Q)^{n-\alpha}} \int_Q |f(y)|  \dy, 
\end{align*} 
where the supremum is taken over all open cubes $Q\subset\Rn$  containing the point $x$ and $\ell(Q)$ denotes the length of the side of the cube. In this article, we consider only open cubes whose sides are parallel to the coordinate axes. When $\alpha = 0$, the operator $\mathcal{M}_0$ coincides with the classical uncentered Hardy–Littlewood maximal function on $\mathbb{R}^n$, which we denote simply by $\mathcal{M}$.

It is well known that $\mathcal{M}$ is bounded on $L^p(\mathbb{R}^n)$ for all $1 < p \leq \infty$, and that it is of weak type $(1,1)$, i.e., bounded from $L^1(\mathbb{R}^n)$ to weak-$L^1(\mathbb{R}^n)$. For $\alpha \in (0,n)$, the fractional maximal function $\mathcal{M}_\alpha$ is bounded from $L^p(\mathbb{R}^n)$ to $L^q(\mathbb{R}^n)$ whenever the pair $(p,q)$ satisfies the relation $\frac{1}{p}=\frac{1}{q}+\frac{\alpha}{n}$, $1<p<\frac{n}{\alpha}$ and $\fma$ maps $L^1(\Rn)$ into the weak-$L^{\frac{n}{n-\alpha}}(\Rn)$ space. Beyond $L^p$-spaces, the behavior of the operators $\mathcal{M}_\alpha$ has also been extensively investigated in other function space settings, including Sobolev spaces  \cites{KS-2003, Weight}, Orlicz spaces \cite{Musil}, and various capacitary function spaces  \cites{ov, tanaka}.
 
This article is devoted to the study of the action of the uncentered fractional maximal function $\fma$ on various spaces of mean oscillation. Furthermore, we consider analogous problems in the framework of dyadic Hausdorff content. John and Nirenberg first explored the concept of functions with bounded mean oscillation (BMO) in their pioneering work \cite{John-Nirenberg}, where they proved the now-classic John–Nirenberg inequality. This result demonstrates that, although such functions may not be bounded, they satisfy an exponential integrability condition. For an open cube $Q_0\subset \mathbb{R}^n$, the space $\cbmo(Q_0)$ is defined to be the space of all locally integrable functions $f$ such that 
\begin{align*}
  \| f \|_{\cbmo(Q_0)} := \sup_{Q} \inf_{c \in \mathbb{R}} \left( \frac{1}{\ell(Q)^n} \int_Q |f(y) - c| \dy \right)<\infty,  
\end{align*}
where the supremum is taken over any finite subcubes $Q$ of $Q_0$ and parallel to $Q_0$.
Originally motivated by problems in elasticity, the BMO space has had a profound impact on harmonic analysis and its applications. BMO includes bounded functions, convolutions of $L^1(\mathbb{R}^n)$ functions with $\log|\cdot|$, and elements of critical Sobolev spaces (see \cite{MS}). Its utility in partial differential equations (PDEs) analysis lies in its ability to capture oscillations that are not adequately described by $L^p$ norms, as illustrated by Brezis–Nirenberg \cite{bn} in their study of elliptic equations with critical exponents. While weak solutions to elliptic or parabolic PDEs may lack $ L^\infty$ bounded gradients or second derivatives, these quantities often exhibit properties that can be described within the BMO framework. This space is characterized by exponential integrability, as established via the John–Nirenberg inequality \cite{John-Nirenberg}. Recent research extends these ideas to lower-dimensional sets, examining the decay of Hausdorff content in their level sets. For a list of related works, we refer to \cites{Adams1973, MS, flz, Cianchi:2008, jud} and the references therein. 

In \cite{Bennet-Sharpley}*{Theorem 4.1}, Bennet-DeVore-Sharpley proved the boundedness of Hardy-Littlewood maximal function $\mathcal{M}$ on $\cbmo(Q_0)$. Later, in \cite{Bennet}*{Lemma 1}, Bennet improved the result of \cite{Bennet-Sharpley}*{Theorem 4.1} by proving that $\mathcal{M}$ maps $\cbmo(Q_0)$ to $\cblo(Q_0)$, where $\cblo(Q_0)$ is the space of bounded lower oscillations and defined by the collection of all locally integrable functions such that
\begin{align*}
  \| f \|_{\cblo(Q_0)} := \sup_{Q} \left( \frac{1}{\ell(Q)^n} \int_Q \left(f(y) - \essinf_{Q}f\right) \dy \right)<\infty,  
\end{align*}
where the supremum is taken over any finite subcubes $Q$ of $Q_0$ and parallel to $Q_0$. In  \cite{fractionalmaximalpaper}*{Corollary 1.2}, Gibara-Kline extended the result on the boundedness of $\cbmo(Q_0)$ to the uncentered fractional maximal function $\mathcal{M}_{\alpha}$ with $0 \leq  \alpha<n$ and proved the following:
\begin{theorem}[Kline-Gibara]\label{Kline-Gibara-BMO-result}
  Let $Q_0\subset \mathbb{R}^n$ be a finite cube and  $0\leq \alpha <n$. If $f \in \cbmo(Q_0)$ with $\supp(f) \subset Q_0$, then $\mathcal{M}_{\alpha}f$ is also in $\cbmo(Q_0)$ and there exist a constant $C>0$ depending only on $\alpha$ and $n$ such that
  \begin{align*}
 \| \fma f \|_{\cblo(Q_0)} \leq C \ell(Q_0)^{\alpha} \| f \|_{\cbmo(Q_0)}. 
 \end{align*}
\end{theorem}
The authors also proved the Theorem \ref{Kline-Gibara-BMO-result} for more general bounded doubling metric measure spaces (see \cite{fractionalmaximalpaper}*{Theorem 1.1}). For recent results related to the boundedness of maximal functions on BMO spaces, we refer to \cites{kara, hlw}.

In this paper, we further improve the boundedness of $\mathcal{M}_{\alpha}$ for $0<\alpha<n$ on $\cbmo(Q_0)$ by showing that it maps $\cbmo(Q_0)$ to $\blo(Q_0)$ with $0<\beta< n$, where  $\blo(Q_0)$ is the space of bounded lower oscillations associated with the dyadic Hausdorff content $\cont$, which we will define in the next few paragraphs.  

Let $\mathcal{D}_0$ be the collection of standard dyadic cubes i.e. $\mathcal{D}_0= \{2^k(l+[0,1)^n): k\in \mathbb{Z}, l \in \mathbb{Z}^n \}$.
For $\beta\in (0,n]$, we define the dyadic Hausdorff content $\cont$ of $E\subset \Rn$  as
\begin{align*}
\cont(E):= \inf\left\{\sum_{i=1}^{\infty}\ell(Q_i)^\beta\, : \, E\subset \cup_{i=1}^\infty Q_i, Q_i\in \mathcal{D}_0 \right\}.
\end{align*} 

Let $\Omega\subset\Rn$ and $f:\Omega  \rightarrow [0,\infty]$. We recall that the Choquet integral with respect to $\cont$ is defined as
		\begin{align}\label{choq-int}
			\int_{\Omega} f \dc := \int_{0}^\infty \cont \left(\{x\in \Omega\,:\, f(x)>t\}\right)\dt,
		\end{align}
where the right-hand side is interpreted as a Lebesgue or improper Riemann integral of the monotone increasing set function
	\begin{align*}
		t \mapsto \cont \left(\{f>t\}\right).
	\end{align*}
   
For open sets $\Omega\subset\Rn$, we define the capacitary integration space  $L^1(\cont, \Omega)$ as the set of all $\cont$-quasicontinuous functions $f$ on $\Omega$ (see p.~6 for the definition) for which 
\begin{align*}
    \|f\|_{L^1(\cont,\Omega)}\vcentcolon= \int_{\Omega} |f|\dc<\infty.
\end{align*}

Let $Q_0 \subset \Rn$ be an open cube and $0<\beta \leq n$. Then a function $f \in L^1(\cont, Q_0)$ is said to be in $\bmo(Q_0)$, the space of a bounded $\beta$-dimensional mean oscillation on $Q_0$, if
\begin{align}\label{beta-bmo}
  \| f \|_{\bmo(Q_0)} := \sup_{Q} \inf_{c \in \mathbb{R}} \left( \frac{1}{\cont(Q)} \int_Q |f(y) - c| \dc(y) \right)<\infty,  
\end{align}
where the supremum is taken over any finite subcubes $Q$ of $Q_0$ and parallel to $Q_0$. The space $\bmo(Q_0)$ was first defined in the paper \cite{Chen-Spector}. For $\beta=n$, the space $\bmob{n}(Q_0)$  coincides with the classical $\cbmo(Q_0)$ space and for  $0<\beta<n$, and $\bmo(Q_0)$ is contained in the space $\cbmo(Q_0)$ (see Lemma \ref{nesting}). For more details about capacitary integration spaces and $\bmo(Q_0)$, we refer to Section \ref{sec:prelim}. 

Let \(Q_0\subset\mathbb{R}^n\) be an open cube and fix \(0<\beta\le n\). For any function $f$, and $A\subset \Rn$, we define 
\begin{align*}
    \enf_A f = \sup\{c\in \mathbb{R}\, :\, \cont(\{x\in A\, :\, f(x)< c\})=0\}.
\end{align*}
When $\beta = n$, this quantity coincides with the essential infimum of $f$ with respect to the Lebesgue measure.
We define the space of \(\beta\)-dimensional bounded lower oscillation on \(Q_0\),
denoted by \(\blo(Q_0)\), as the set of all real-valued functions
\(f\in L^{1}(\cont,Q_0)\) for which  
\begin{align}\label{blo-def}
\|f\|_{\blo(Q_0)}=\sup_{Q}\frac{1}{\cont(Q)}\int_{Q}\left(f(x)-\enf_{y \in Q} f(y)\right)\dc(x) <\infty,
\end{align}
where the supremum is taken over any finite subcubes $Q$ of $Q_0$ and parallel to $Q_0$. Clearly, $\blo(Q_0) \subset \bmo(Q_0)$ for $0<\beta\leq n$ and for $\beta=n$, the space $\blob{n}(Q_0)$  coincides with the classical $\cblo(Q_0)$ space. In Section \ref{Sec:BLO-JN}, we establish that functions in the space $\blo(Q_0)$ with $0<\beta\leq n$, satisfy a John–Nirenberg type inequality with respect to the $\blo$ norm introduced in  \eqref{blo-def}. This inequality plays a pivotal role in the proof of Theorem \ref{MainThm2} and is also of independent interest in the study of function spaces.

The following theorem is the first main result of our paper.
\begin{theorem}\label{MainThm2}
Let $0 < \alpha < n$, $0 < \beta \leq n$, and let $Q_0 \subset \Rn$ be a finite cube. Assume $f \in \cbmo(Q_0)$ with $\supp(f) \subset Q_0$. Then there exists a constant $C > 0$, depending only on $\alpha$, $\beta$, and $n$, such that
\begin{align*}
 \| \fma f \|_{\blo(Q_0)} \leq C \ell(Q_0)^{\alpha} \| f \|_{\cbmo(Q_0)}.   
\end{align*}
\end{theorem}

Note that from Corollary \ref{nesting-blo}, we conclude $\blo(Q_0) \subset \cblo(Q_0)$ for $0<\beta<n$. Therefore, Theorem \ref{MainThm2} is an improvement of Theorem \ref{Kline-Gibara-BMO-result}. Our argument for proving Theorem \ref{MainThm2} is inspired by \cite{fractionalmaximalpaper}*{Theorem 1.1}. To get the improvement in terms of $\beta$-dimensional $\blo$ norm, we incorporated the boundedness result of $\fma$ from \cite{Sawyer}*{Theorem B}. When $0<\alpha<1$ and $f \in \cbmo(Q_0)$, it is known that $\fma f$ is H\"older continuous (see \cite{hln}*{Theorem~3.1}). However, it is not evident how Theorem \ref{MainThm2} could be deduced solely from the H\"older continuity of the fractional maximal function over the full range $0<\alpha<n$. 

It is worth noting that Theorem~\ref{MainThm2} does not cover the case $\alpha = 0$. In fact, we cannot expect Theorem~1.2 to hold for the Hardy–Littlewood maximal function $\mathcal{M}$. For example, let $H_k$ be a $k$-hyperplane in $\mathbb{R}^n$, where $k$ is a natural number with $0 < k < n$, and consider the function $f(x) = \log(\dist(x, H_k))$. One can show that $f \in \cbmo(\Rn)$, but $\mathcal{M}f(x) = +\infty$ for $\mathcal{H}^{k}_{\infty}$ a.e. $x \in H_k$. Hence, $\mathcal{M}f \notin \cbmo^{ k}(\Rn)$. Nevertheless, we can establish the following alternative improvement in our next result.

\begin{theorem}\label{asthe}
    For $ 0 < \beta \leq n$, there exists a constant $C>0$, depending only on $\beta$, and $n$ such that
\begin{align*}
\|\mathcal{M} f\|_{\blob{\beta}(Q)} \leq C \|f\|_{\bmob{\beta}(Q)},
\end{align*}
 where $Q$ is either a finite cube or $\mathbb{R}^n$.
\end{theorem} 

In the above, the spaces $\bmob{\beta}(\mathbb{R}^n)$ and $\blob{\beta}(\mathbb{R}^n)$ are defined analogously to \eqref{beta-bmo} and \eqref{blo-def}, respectively, by taking the supremum over all cubes in $\mathbb{R}^n$. In fact, we can prove a general version of Theorem~\ref{asthe} (see Theorem~\ref{MainThm} in Section~\ref{Sec:beta-BMO-Maximal function}) by establishing the result for the $\beta$-dimensional Hardy–Littlewood maximal function $\mathcal{M}^{\beta}$ (see \eqref{man_def} for the definition), which was introduced in \cite{ChenOoiSpector}.

Another important mean oscillation space, which is a subset of BMO, is the space of vanishing mean oscillation (VMO). This VMO space is introduced by Sarason in \cite{Sarason} as the closure of $C_c^\infty(\Rn)$ in the BMO norm, refine BMO spaces by capturing functions whose mean oscillation vanishes at small scales, thus bridging bounded oscillation and continuity. An equivalent definition of VMO is the following: Let $Q_0\subset \Rn$ be an open cube. Then a function $f\in \cbmo(Q_0)$ is said to be in $\cvmo(Q_0)$ if 
\begin{align*}
    \lim_{r\rightarrow 0^+ } \sup_{Q \subset Q_0\, , \ell(Q)\leq r}	\O(f,Q)=0,
\end{align*}
where, $\O(f,Q)$ is the modulus of mean oscillation of $f$ on $Q$, defined by
\begin{align*}
 \O(f,Q)= \inf_{c\in\R} \frac{1}{\ell(Q)^n}\int_{Q}|f-c|\dx.   
\end{align*}

VMO functions, while not necessarily continuous, exhibit continuity-like behavior locally, making them crucial in harmonic analysis, PDEs, and functional analysis, particularly when fine regularity or compactness is needed. For instance, solutions to elliptic or parabolic PDEs in divergence form can remain well-behaved even with coefficients lacking full continuity, as shown in \cite{cfl-2}. Calderón–Zygmund operators preserve VMO, a result due to Coifman-Rochberg-Weiss \cite{crw}, which underpins much of their utility in singular integral theory. Recent developments, such as those by Nowak \cite{noak}, address higher integrability and differentiability in nonlocal equations with irregular, VMO-type, or small BMO coefficients. Similarly, Hofmann-Mayboroda \cite{hm} explored VMO within Hardy spaces adapted to elliptic operators with rough coefficients, reflecting its modern significance in harmonic analysis. In this work, we define VMO using Hausdorff content, namely $\vmo$, and this opens up a new exploration of the above-mentioned studies to irregular domains and fractal settings.

Now, we state our third main result, which constitutes an improvement of the work of Gibara–Kline \cite{fractionalmaximalpaper}*{Corollary~1.2}. Let $Q_0$ be a finite cube in $\Rn$, and let $0<\alpha<n$. It was shown in \cite{fractionalmaximalpaper} that if $f\in \cvmo(Q_0)$, then $\fma f\in \cvmo(Q_0)$. In the following, we extend this conclusion by demonstrating that it remains valid under the broader assumption that $f\in \cbmo(Q_0)$, where $\cbmo(Q_0)$ properly contains $\cvmo(Q_0)$. The precise statement is presented below.	
\begin{theorem}\label{VMO-improvement}
Let $0 < \alpha < n$ and $Q_0\subset \mathbb{R}^n$ be a finite open cube. Let $f\in \cbmo(Q_0)$ and $\supp(f)\subset Q_0$. Then $\fma f$ is uniformly continuous on $\Rn$. Consequently, 
\begin{align*}
    \fma f\in\cvmo(Q_0).
\end{align*}
\end{theorem}

In our final main result, we investigate the action of $\fma$ on the vanishing $\beta$-dimensional mean oscillation space associated with the dyadic Hausdorff content $\cont$. This space, which captures finer oscillatory behavior relative to the underlying content, will be defined precisely below. Our goal is to establish that $\fma$ preserves the vanishing $\beta$-dimensional mean oscillation property under suitable conditions, thereby extending the regularity theory in this setting.

Let $0< \beta \leq n$ and $Q_0$ be any finite cube of $\Rn$ and take $f\in L^1(\cont, Q_0)$. Then we define the \textit{$\beta$-dimensional mean oscillation} of $f$ on $Q\subset Q_0$ as follows
\begin{align*}
\O_{\beta}(f,Q):=\inf_{c \in \mathbb{R}}\frac{1}{\cont(Q)}\int_{Q}|f-c|\dc,
\end{align*}
and for $r>0$, the $\beta$-dimensional modulus of oscillation of $f$ as follows
\begin{align*}
\omega_{\beta}(f, r):= \sup_{Q \subset Q_0\,,  \ell(Q)\leq r}	\O_{\beta}(f,Q).
\end{align*}
Then we say that $f\in \bmo(Q_0)$ has \textit{vanishing $\beta$-dimensional mean oscillation} on $Q_0$ (denoted as $f\in \vmo(Q_0)$)  if 
	\begin{align*}
		\omega_{\beta}(f,0):=\lim_{r\rightarrow 0^+ } \omega_{\beta}(f,r) =0.
	\end{align*}
Note that for $\beta=n$, $\cvmo^n(Q_0)$ coincides with $\cvmo(Q_0)$ and  $\cvmo^{\beta}(Q_0)\subset \cvmo(Q_0)$ for $0<\beta <n$ (see Lemma \ref{vmo-comp}).
 
In the following theorem, we address the question of whether the vanishing $\beta$-dimensional mean oscillation property is preserved under the action of uncentered fractional maximal functions. Let $Q_0 \subset \Rn$ be a finite cube, and suppose that $f \in \vmo(Q_0)$ with $\fma f$ is not identically equal to infinity. Since $\bmo(Q_0) \subset \cbmo(Q_0)$ for $0 < \beta \leq  n$ (see Lemma~\ref{nesting}), Theorem~\ref{VMO-improvement} implies that $\fma f \in \cvmo(Q_0)$. This naturally leads to the question of whether $\fma f \in \vmo(Q_0)$ for $0<\beta<n$. We answer this question affirmatively, thereby establishing the stability of the vanishing $\beta$-dimensional mean oscillation property under the action of the uncentered fractional maximal functions. 
\begin{theorem}\label{thm-VMO-main}
Let $0 < \beta  \leq   n$, $\alpha \in [0,\beta)$ and $Q_0 \subset \mathbb{R}^n$ be a finite open cube. Let $f\in \vmo(Q_0)$ and $\supp f\subset Q_0$. Then if $\mathcal{M}_\alpha f$ is not identically infinite, then $\mathcal{M}_\alpha f\in \vmo(Q_0)$.	
\end{theorem}
\begin{remark}
  In \cite{xianbao}*{Theorem 3.1}, Shaabani established that the uncentered Hardy–Littlewood maximal function $\mathcal{M}$ preserves the space $\cvmo(\mathbb{R}^n)$. By careful observation of the proof of Theorem \ref{thm-VMO-main}, we observed that for $\alpha=0$, we can extend the mapping properties of $\mathcal{M}$ on $\vmo(\Rn)$.
\end{remark}
    
The structure of the remainder of the paper is as follows. In Section~\ref{sec:prelim}, we recall essential preliminaries concerning dyadic Hausdorff content and Choquet integrals, and establish some basic results that will be used in the proofs of the main theorems. Section~\ref{Sec:BLO-JN} is devoted to the proof of a John–Nirenberg-type inequality for the space $\blo$, along with several of its consequences. In Section~\ref{Sec:BMO-improve}, we establish our first main result, Theorem~\ref{MainThm2}. Section~\ref{Sec:beta-BMO-Maximal function} focuses on the boundedness of the $\beta$-dimensional uncentered maximal function on $\bmo$, concludes with the proof of Theorem~\ref{MainThm}. In Section~\ref{Sec: VMO-boundedness}, we prove Theorem~\ref{VMO-improvement}, and in Section~\ref{Sec:Beta-VMO}, we complete the proof of Theorem~\ref{thm-VMO-main}, which concerns the mapping properties of the uncentered fractional maximal function on $\vmo$. Finally, in Appendix~\ref{appen}, we prove results related to quasicontinuity of maximal functions.

\textbf{Notations:} In this article, we use the following notations:
\begin{itemize}
    \item We use $a\lesssim b$ to mean that there exists a positive universal constant $C>0$ such that $a\leq C b$. Additionally, $a\cong b$ indicates that $a\lesssim b$ and $b\lesssim a$ together. A constant is universal if it depends on the Euclidean dimension $n$ and $\beta$ associated with the Hausdorff content.
    \item In this article, we will always consider open cubes if it is not specified. 
    \item Here we denote cube $Q(x,r)$ as an open cube with center at $x$ and with length $r>0$.
    \item For any set $A\subset \Rn$, we write $|A|$ is the Lebesgue measure of $A$.
    \item For any cube $Q$ on $\Rn$, and any nonnegative function $f$, we write 
\begin{align}\label{averageQrespecttobeta}
  \avgd{Q}:= \frac{1}{\cont(Q)}\int_{Q} f\dc,
\end{align}
and 
\begin{align*}
  \avgl{Q} :=\frac{1}{\ell(Q)^n}\int_{Q}f\dx, 
\end{align*}
where $\dx$ is the volume element associated with the usual Lebesgue measure.
\end{itemize}

\section{Preliminaries} \label{sec:prelim}
Let $\mathcal{P}(\Rn)$ denote the collection of all subsets of $\Rn$. We recall that the dyadic Hausdorff content $\cont:\mathcal{P}(\Rn)\rightarrow [0,\infty]$ is an outer measure satisfying the following properties:
\begin{itemize}
	\item[(i)] (Null set) $\cont(\emptyset)=0$;
	\item[(ii)] (Monotonicity) If $A\subseteq B$, then $\cont(A)\leq \cont(B)$;
	\item[(iii)] (Countable subadditivity) If $A\subseteq \bigcup_{i=1}^\infty A_i$, then $\cont(A)\leq \sum_{i=1}^\infty\cont(A_i)$;
	\item[(iv)] (Strong subadditivity) $\cont(A\cup B) + \cont(A\cap B) \leq \cont(A) + \cont(B)$;
	\item[(v)] (Doubling property) There exists a universal positive constant $C$ such that for every $A\in \mathcal{P}(\Rn)$,
	\begin{align*}
	    \cont(2A)\leq C\cont(A), \quad \text{where}\quad 2A:=\{2a\, :\,a\in A\}.
	\end{align*}
\end{itemize}
We also define the Hausdorff content associated with general cubes (not necessarily dyadic cubes), denoted by $\widetilde{\mathcal{H}}^\beta_\infty$, as follows:
\begin{align*}
    \widetilde{\mathcal{H}}^\beta_\infty (E) := \inf\left\{\sum_{i=1}^{\infty}\ell(Q_i)^\beta : E\subseteq \bigcup_{i=1}^{\infty} Q_i \right\}.
\end{align*}
It follows that both the contents are equivalent, i.e., there exist positive constants $C_1$ and $C_2$, depending only on $\beta$, such that
\begin{align}\label{equiv-cont}
C_1\widetilde{\mathcal{H}}^\beta_\infty(E)\leq \mathcal{H}_\infty^\beta(E)\leq C_2\widetilde{\mathcal{H}}^\beta_\infty(E)
\end{align}
(see \cite{Yang-Yuan}*{Proposition~2.3}). Since $\widetilde{\mathcal{H}}^\beta_\infty (Q) = \ell(Q)^\beta $ for any cube $Q \subset \mathbb{R}^n$, inequality \eqref{equiv-cont} immediately yields
\begin{align}\label{equiv-cont-cube}
C_1\ell(Q)^\beta\leq \mathcal{H}_\infty^\beta(Q)\leq C_2\ell(Q)^\beta.
\end{align}

We next present several fundamental properties of the Choquet integral with respect to the dyadic Hausdorff content defined in \eqref{choq-int}. The proofs of the following lemma follow directly from the arguments presented in \cite{AdamsChoquet} (see also \cite{CMS} and \cite{Herjulehto-Petteri_2024}*{Lemma~3.2}), combined with the strong subadditivity of the dyadic Hausdorff content $\mathcal{H}_\infty^\beta$.
\begin{lemma}\label{basicChoquetintegral}
Suppose $\Omega$ is an open subset of $\Rn$ and $\beta\in (0,n]$. Then the following properties hold:
\begin{enumerate}
\item If $c \geq 0$ and $f \geq 0$, then
\begin{align*}
\int_\Omega c f\dc  = c \int_\Omega f\dc.
\end{align*}

\item \label{sub} For nonnegative functions $f_n$, we have
\begin{align*}
\int_\Omega \sum_{n=1}^\infty f_n\dc \leq \sum_{n=1}^\infty \int_\Omega f_n\dc.
\end{align*}

\item Let $1 \leq p < \infty$ with $\frac{1}{p}+\frac{1}{p'}=1$. Then for nonnegative functions $f$ and $g$, we have
\begin{align*}
\int_\Omega f(x)g(x)\dc \leq \left(\int_\Omega |f(x)|^p\dc\right)^{\frac{1}{p}} \left(\int_\Omega |g(x)|^{p'}\dc\right)^{\frac{1}{p'}}.
\end{align*}
\end{enumerate}
\end{lemma}

Let $\Omega\subseteq \Rn$ be an open set and $f:\Omega\rightarrow\R$ be a function. We say that a function $f$ is $\cont$-quasicontinuous on $\Omega$ if for every $\epsilon>0$, there exists an open set $O_\epsilon\subset\Omega$ such that $\cont(O_\epsilon)<\epsilon$ and $f|_{\left(\Omega\setminus O_\epsilon\right)}$ is continuous. For $\Omega\subset\Rn$ and $1\leq p<\infty$, we define
    \begin{align}\label{l1}
        L^p(\cont, \Omega) \vcentcolon= \left\{ f\,: f \text{ is } \cont\text{-quasicontinuous on }\Omega \text{ and } \|f\|_{L^p(\cont,\Omega)}\vcentcolon= \left( \int_{\Omega} |f|^p\dc \right)^{\frac{1}{p}}<\infty \right\}.
    \end{align}
The space $L^p(\cont, \Omega)$ is a Banach space, and the intersection $L^p(\cont, \Omega) \cap C_b(\Omega)$, where $C_b(\Omega)$ denotes the space of continuous and bounded functions on $\Omega$, is dense in $L^p(\cont, \Omega)$. Also for $\beta=n$, we write $L^p(\contb{n}, \Omega)=L^p(\Omega)$, the set of usual $p$ integrable functions.  

Let $Q_0$ be an open cube in $\Rn$. For $0<\beta \leq n$, $1\leq p<\infty$, a function $f $ is said to be in $\bmobp{\beta}{p}(Q_0)$ if $f$ is $\cont$-quasicontinuous on $Q_0$ and
\begin{align}\label{bmo-1}
  \| f \|_{\bmobp{\beta}{p}(Q_0)} := \sup_{Q \subseteq Q_0} \inf_{c \in \mathbb{R}} \left( \frac{1}{\cont(Q)} \int_Q |f(y) - c|^p \dc(y) \right)^{1/p}<\infty.  
\end{align}
The analogue of the John–Nirenberg inequality for \(\cbmo^{\beta}(Q_{0})\) proved in
\cite{Chen-Spector}*{Theorem~1.3} yields the norm equivalence
\begin{equation}\label{equi-p-bmo}
    \frac{1}{C}\,\|f\|_{\cbmo(Q_0)}
    \;\le\;
    \|f\|_{\cbmo^{\beta,p}(Q_0)}
    \;\le\;
    C\,p\,\|f\|_{\cbmo(Q_0)},
\end{equation}
for every \(p \geq 1\) and all \(f \in \cbmo^\beta(Q_0)\) (see \cite{Chen-Spector}*{Corollary 1.5}).

We next record the following lemma, which can be found in \cite{Chen-Spector}*{Corollary 1.6}.

\begin{lemma}\label{nesting}
Let $0<\gamma \leq \beta \leq n$ and suppose $u \in \bmob{\gamma}(Q)$, where $Q$ is finite cube or $\Rn$. Then $u \in \bmo(Q)$ and
\begin{align*}
\|u\|_{\bmo(Q)} \leq C\|u\|_{\bmob{\gamma}(Q)}
\end{align*}
for some positive constant $C=C(\gamma,\beta,n)$ independent of $u$.
\end{lemma}

The next two results, whose proofs are essentially identical and rely on the non‑trivial sublinearity of the Choquet integral with respect to the dyadic Hausdorff content, will be used repeatedly in the subsequent sections.

\begin{lemma}\label{tri_cont}
Let \(0 < \beta \leq n\), and let \(Q\) be a cube in \(\mathbb{R}^n\). Then, for any \(c \geq 0\) and nonnegative function $f$, one has
\begin{align*}
    | f_{Q,\beta} - c | \leq \frac{1}{\cont(Q)} \int_Q |f - c| \dc,
\end{align*}
where \(f_{Q,\beta}\) is defined as in \eqref{averageQrespecttobeta}.
\end{lemma}

	\begin{proof}
		By Lemma~\ref{basicChoquetintegral}, we have
        \begin{align}\label{averageestimatefirstdirection}
      \nonumber  f_{Q,\beta} &= \frac{1}{\cont(Q)} \int_{Q}|f-c +c |\dc \\
       \nonumber     &\leq 	 \frac{1}{\cont(Q)}\int_{Q}|f-c|\dc+  \frac{1}{\cont(Q)}\int_{Q}c\dc\\
            & =  \frac{1}{\cont(Q)}\int_{Q}|f-c|\dc+ c .
        \end{align}             
		Now subtracting both sides by $c$, we get
\begin{align*}
			f_{Q,\beta} -c \leq  \frac{1}{\cont(Q)}\int_{Q}|f-c|\dc.
		\end{align*}
A similar argument as in~\eqref{averageestimatefirstdirection} yields
\begin{align*}
    c-f_{Q,\beta}\leq  \frac{1}{\cont(Q)} \int_{Q}|f-c|\dc.
\end{align*}
The combination of the last two inequalities then gives the desired estimate.
\end{proof}

As a consequence of Lemma~\ref{tri_cont}, we have the following equivalence of $\bmo$ norms for nonnegative functions.
\begin{theorem}\label{bmo-remark-def}
  For $0<\beta\leq n$ and any nonnegative function $f\in \bmo(Q_0)$, we have
\begin{align*}
  \| f \|_{\bmo(Q_0)} \cong \sup_{Q}  \left( \frac{1}{\cont(Q)} \int_Q |f(y) - \avgd{Q}| \dc(y) \right),  
\end{align*}
where the supremum is taken over any finite subcubes $Q$ of $Q_0$ and parallel to $Q_0$. 
\end{theorem}

Let us define the $\beta$-dimensional maximal function $\mathcal{M}^{\beta}$ associated with $\cont$. For $0<\beta \leq n$ and $f\in L^1_{loc}(\cont, \Rn)$, we define 
\begin{align}\label{man_def}
\mathcal{M}^{\beta} f (x):=\sup_{Q \ni x} \frac{1}{\cont(Q)} \int_{Q}|f(y)|\dc(y). 
\end{align} 
The following lemma captures a natural continuity property of the $\beta$-dimensional maximal function with respect to pointwise differences
\begin{lemma}
    Let $0  < \beta \leq n$ and $f,$ $g $ be functions of $\mathbb{R}^n$. If $\mathcal{M}^\beta f (x)<\infty$ and $\mathcal{M}^\beta g (x)<\infty$ for some $x \in \mathbb{R}^n,$ then
    \begin{align}\label{diff_max}
        \left| \mathcal{M}^\beta f (x) - \mathcal{M}^\beta g (x)  \right| \leq \mathcal{M}^\beta (f-g) (x ).
    \end{align}
\end{lemma}
 \begin{proof}
     Without loss of generality, assume that \(\mathcal{M}^\beta f(x) \ge \mathcal{M}^\beta g(x)\). By Lemma \ref{basicChoquetintegral} and the definition of $\beta$-dimensional maximal function associated to $\cont$, we have for every cube $Q$ containing $x$, 
\begin{align*}
    \frac{1}{ \cont (Q)} \int_Q |f (y)| \dc (y) & \leq \frac{1}{ \cont (Q)} \int_Q |f (y) - g(y)| \dc (y ) + \frac{1}{ \cont (Q)} \int_Q | g(y)| \dc (y ) \\
    & \leq \mathcal{M}^\beta (f-g) (x) +\mathcal{M}^\beta g (x).
\end{align*}
Taking the supremum over all such cubes \(Q\) yields
\begin{align*}
   \mathcal{M}^\beta f (x) \leq \mathcal{M}^\beta (f-g) (x) + \mathcal{M}^\beta g (x).
\end{align*}
The proof is completed by interchanging the roles of $f$ and $g$.
\end{proof}

\begin{lemma}\label{abs-bmo}
Let \(0 < \beta \le n\) and let \(Q \subset  \mathbb{R}^{n}\) be a finite cube or $\Rn$.
If \(f \in \bmo(Q)\), then \(|f| \in \bmo(Q)\).
Moreover,
\[
   \bigl\|\,|f|\,\bigr\|_{\bmo(Q)}
   \;\le\;
   2\,\|f\|_{\bmo(Q)}. 
\]
\end{lemma}

\begin{proof}
		Fix a cube \(Q' \subseteq Q\) and choose
		\begin{align*}
			c_{Q'}= \argmin_{c \in \mathbb{R}} \frac{1}{\cont(Q')}\int_{Q'} \, |f-c| \dc.
		\end{align*}
	By definition of the \(\cbmo^{\beta}\)-norm, we have \begin{align}\label{argbmobeta}
	    \frac{1}{\cont(Q')}\int_{Q'}|f-c_{Q'}| \dc \leq \|f\|_{\bmo}.
	\end{align}
		For every \(x \in Q'\), Lemma~\ref{tri_cont} and (\ref{argbmobeta}) yield that
		\begin{align*}
			& \left| |f(x)|-\frac{1}{\cont(Q')}\int_{Q'}|f(y)| \dc (y) \right| \\
			\leq & \frac{1}{\cont(Q')}\int_{Q'} |f(x)-f(y)| \dc(y) \\
			=& \frac{1}{\cont(Q')} \int_{Q'} |f(x)-c_{Q'}+c_{Q'}-f(y)| \dc (y) \\
			\leq & |f(x)-c_{Q'}| + \frac{1}{\cont(Q')} \int_{Q'} |f(y)-c_{Q'}| \dc(y)\\
			\leq & |f(x)-c_{Q'}| + \|f\|_{\bmo(Q)}.
		\end{align*}
			Integrating the preceding estimate over \(Q'\)  with respect to the dyadic Hausdorff content \(\cont\) and dividing by \(\cont(Q')\), we obtain
		\begin{align*}
			&\frac{1}{\cont(Q')} \int_{Q'} \left| |f(x)|-\frac{1}{\cont(Q')}\int_{Q'}|f| \dc \right| \dc (x) \\ \leq &  \frac{1}{\cont(Q')} \int_{Q'} |f(x)-c_{Q'}| \dc (x) + \|f\|_{\bmo(Q)}\\
			\leq  & 2 \|f\|_{\bmo(Q)}.
		\end{align*}
		
		Therefore, for each cube $Q'\subset Q$, we have
		\begin{align}\label{eq-infy-bmo}
			&\inf_{c\in \mathbb{R}} \frac{1}{\cont(Q')} \int_{Q'} \left| |f(x)|-c\right| \dc (x) \\ \leq& \frac{1}{\cont(Q')} \int_{Q'} \left| |f(x)|-\frac{1}{\cont(Q')}\int_{Q'}|f| \dc \right| \dc (x)\nonumber\\ \leq &
			 2 \|f\|_{\bmo(Q)}.   \nonumber
		\end{align}
		Taking the supremum over all such cubes gives
		\begin{align*}
			\bigl\|\,|f|\,\bigr\|_{\cbmo^{\beta}(Q)}\leq 2 \|f\|_{\bmo(Q)},
		\end{align*}
which proves the lemma.
	\end{proof}

Fix a cube \(Q_{0}\subset\mathbb{R}^{n}\), \(\beta\in(0,n]\) and \(1\le p<\infty\).
For \(x\in Q_{0}\) and \(r>0\) denote by \(Q(x,r)\) the cube in \(\mathbb{R}^{n}\) centered at \(x\) with side‑length \(r\).
We define \(\overline{\cbmo}^{\,\beta,p}(Q_{0})\) as the collection of all \(\cont\)-quasicontinuous functions \(f\) on \(Q_{0}\) for which
\begin{equation}\label{bmo-2}
  \|f\|_{\overline{\cbmo}^{\,\beta,p}(Q_{0})}
  \vcentcolon=
  \sup_{x\in Q_{0},\,r>0}
  \inf_{c\in\mathbb{R}}
  \biggl(
      \frac{1}{\cont(Q(x,r))}
      \int_{Q(x,r)\cap Q_{0}}
            |f(y)-c|^{p}\dc(y)
  \biggr)^{1/p}
  <\infty.
\end{equation}
 In a similar spirit, we introduce the \emph{$\beta$-dimensional bounded lower oscillation} space \(\cblo^{\beta,p}(Q_{0})\) as the set of all \(\cont\)-quasicontinuous functions \(f\) on \(Q_{0}\) for which
	\begin{align}\label{blo}
\|f\|_{\cblo^{\beta,p}(Q_0)}=\sup_{Q^\prime}\left(\frac{1}{\cont(Q^\prime)}\int_{Q^\prime}\left(f(x)-\enf_{y \in Q^\prime} f(y)\right)^p\dc(x) \right)^{\frac{1}{p}}<\infty,
	\end{align}
	where the supremum is taken over any finite subcubes $Q^\prime$ of $Q_0$. Moreover, for $p=1$, we denote $\cblo^{\beta,1}=\cblo^{\beta}$. Finally, define \(\overline{\cblo}^{\,\beta}(Q_{0})\) as the collection of all locally integrable
\(f\) satisfying
    \begin{align}\label{blo-2}
     \|f\|_{\overline{\cblo}^\beta(Q_0)}= \sup_{Q(x,r), x \in Q_0} \frac{1}{\cont(Q)}\int_{Q(x,r)\cap Q_0}\left(f(z)-\enf_{y \in Q^\prime} f(y)\right)\dc(z) <\infty,
    \end{align}    
    and this definition is equivalent with $\blo(Q_0)$ (see \eqref{ywl-blo}). Moreover, we will prove in Remark~\ref{equi-bmo} that
\[
   \overline{\cbmo}^{\,\beta,p}(Q_{0}) = \cbmo^{\beta,p}(Q_{0}),
 \text{ and }\, 
   \overline{\cblo}^{\,\beta}(Q_{0}) = \cblo^{\beta}(Q_{0}).
\]
Thus, the local and global definitions coincide in both the \(\cbmo^{\beta,p}\) and \(\cblo^{\beta}\) settings.

Next, we record a result that captures a fundamental point in the analysis of the Choquet integral with respect to the dyadic Hausdorff content. This result appears in \cite{ov}, and we also refer the reader to \cite{BCRS}*{Proposition 2.8} for a proof.
	\begin{proposition}\label{deco-cube-int}
Let $\{Q_{j}\}_j$ be a non-overlapping collection of dyadic cubes in $\mathcal{D}_0$. Suppose there exists a constant $A \geq 1$ such that
\begin{equation}\label{pack-int-hyp-1}
			\sum_{Q_{j} \subset Q} \cont\left(Q_{j}\right) \leq A \, \cont(Q) 
		\end{equation}	
        holds for each $Q \in \mathcal{D}_0$.
Then, for any nonnegative function $f$, it holds that
		\begin{equation}\label{pack-int-hyp-2}
			\sum_{j} \int_{Q_{j}} f \dc \leq A \int_{\cup_{j} Q_{j}} f \dc.
		\end{equation}
\end{proposition}

\begin{lemma}\label{diffBMO}
Let $\beta \in (0, n]$, $\gamma \in [0,\beta)$, $p \in [1, \infty)$, and let $f \in L^p(\cont,Q_0)$ be a nonnegative function with $\mathrm{supp}(f) \subset Q_0$. Then for every cube $Q$ satisfying $Q \cap Q_0 \neq  \emptyset $, there exists a cube $P_Q$ with the following properties:
\begin{enumerate}[label=(\roman*)]
\item $P_Q \subseteq Q_0 \cap 2Q$;  
\item $P_Q \supseteq Q_0 \cap Q$;  
\item The side length of $P_Q$ satisfies
   \begin{align*}
   \ell(P_Q) = \begin{cases}
   \ell(Q), & \text{if } \ell(Q) < \ell(Q_0),\\[6pt]
   \ell(Q_0), & \text{if } \ell(Q) \geq \ell(Q_0);
   \end{cases}
   \end{align*}
\item For every $c\in\R$, there exists a constant $D$ depending only on $\beta$, $\gamma$ and $p$ such that 
   \begin{align*}
   \left( \frac{1}{\cont(Q)^{1-\frac{\gamma}{\beta}}} \int_{Q \cap Q_0} |f - c|^p \dc \right)^{\frac{1}{p}}
   \leq  D \left( \frac{1}{\cont(P_Q)^{1-\frac{\gamma}{\beta}}} \int_{P_Q} |f - c|^p \dc \right)^{\frac{1}{p}}.
   \end{align*}
   Moreover, when \(\beta = n\), the above inequality holds with \(D = 1\).
\end{enumerate}

\end{lemma}
\begin{proof}
We consider two cases based on the relative sizes of $Q_0$ and $Q$.

\textbf{Case 1:} $\ell(Q)<\ell(Q_0)$.  
Since $Q\cap Q_0\neq \emptyset$, we can choose a cube $P_Q \subseteq Q_0$ with $\ell(P_Q) = \ell(Q)$ such that $P_Q \supseteq Q \cap Q_0$. With this choice of $P_Q$, we have
\begin{align*}
\left( \frac{1}{\cont(Q)^{1-\frac{\gamma}{\beta}}} \int_{Q \cap Q_0} |f - c|^p \dc \right)^{\frac{1}{p}}
&\leq \left( \frac{1}{\cont(Q)^{1-\frac{\gamma}{\beta}}} \int_{P_Q} |f -c|^p \dc \right)^{\frac{1}{p}}\\
& = \left( \frac{\cont(P_Q)^{1- \frac{\gamma}{\beta}}}   {\cont(Q)^{1- \frac{\gamma}{\beta}}} \right)^{\frac{1}{p}} \left( \frac{1}{\cont(P_Q)^{1-\frac{\gamma}{\beta}}} \int_{P_Q} |f -c|^p \dc \right)^{\frac{1}{p}}\\
&\leq \left( \frac{C_2}{C_1} \right)^{(1-\gamma/\beta)/p} \left( \frac{1}{\cont(P_Q)^{1-\frac{\gamma}{\beta}}} \int_{P_Q} |f -c|^p \dc \right)^{\frac{1}{p}},
\end{align*}
where \(C_{1},C_{2}>0\) are the constants depending only on \(\beta\) that appear in~\eqref{equiv-cont-cube}. Hence \(\ell(P_{Q})=\ell(Q)\) and the inclusion
\(P_{Q}\subset Q_{0}\cap 2Q\) verify all the required properties in this case.

\textbf{Case 2:} $\ell(Q)\geq \ell(Q_0)$.  
In this case, a similar argument as given in Case 1 gives
\begin{align*}
\left( \frac{1}{\cont(Q)^{1-\frac{\gamma}{\beta}}} \int_{Q \cap Q_0} |f - c|^p \dc \right)^{\frac{1}{p}}
& \leq  \left( \frac{\cont(Q_0)^{1- \frac{\gamma}{\beta}}}   {\cont(Q)^{1- \frac{\gamma}{\beta}}} \right)^{\frac{1}{p}} \left( \frac{1}{\cont(Q_0)^{1-\frac{\gamma}{\beta}}} \int_{Q_0} |f -c|^p \dc \right)^{\frac{1}{p}}\\
&\leq \left( \frac{C_2}{C_1} \right)^{(1-\gamma/\beta)/p} \left( \frac{1}{\cont(Q_0)^{1-\frac{\gamma}{\beta}}} \int_{Q_0} |f -c|^p \dc \right)^{\frac{1}{p}}.
\end{align*}
  Hence, the choice \(P_{Q}:=Q_{0}\) satisfies all the required properties and yields the desired estimate.

If \(\beta=n\), then \(\cont(Q') = \ell(Q')^{\,n}\) for every cube $Q'$ of $\mathbb{R}^{n}$. Consequently, one may take \(C_{1}=C_{2}=1\) in~\eqref{equiv-cont-cube}, and thus the inequality holds with \(D=1\), completing the proof. 
\end{proof}
 This suggests that for \(\beta=n\) (so that \(D=1\)), we have the formula
\begin{align}\label{ywl-mf}
  \mathcal{M}_\alpha  f(x) 	= \sup_{\{ Q' \text{ is subcube of }Q_0, \,x\in Q^\prime\}} \frac{1}{\ell(Q')^{n-\alpha }}\int_{Q^\prime}|f(y)| \dy  
\end{align}
holds for any $x \in Q_0$, whenever $\supp(f)\subset Q_0$.
\begin{remark}\label{equi-bmo}
We can immediately observe the following facts from Lemma~\ref{diffBMO}. 
\begin{enumerate}
    \item[(i)] For any $\beta \in (0, n]$, and $1\leq p <\infty$, we have
    \begin{align}\label{ywl-bmo} 
        \|f \|_{\cbmo^{\beta,p}(Q_0)} \cong \| f \|_{\overline{\cbmo}^{\beta,p}(Q_0)}.
    \end{align}
    \item[(ii)] For any $\beta \in (0, n]$, we have
    \begin{align}\label{ywl-blo} 
        \|f \|_{\blo(Q_0)} \cong \| f \|_{\overline{\cblo}^{\beta}(Q_0)}.
    \end{align}
\end{enumerate}    
\end{remark}
 \section{John-Nirenberg inequality for \texorpdfstring{$\beta$}{beta}-dimensional bounded lower oscillation}\label{Sec:BLO-JN}
 This section is devoted to the John–Nirenberg inequality for the space $\blo(Q_0)$. Our result is inspired by the work of \cite{WangZhouTeng}*{Theorem~2.2}. However, the techniques employed in \cite{WangZhouTeng} are not directly applicable in our setting due to the nonlinearity of the Choquet integral with respect to dyadic Hausdorff content $\cont$. To address this issue, we adopt the approach developed in \cite{Chen-Spector} to establish the John-Nirenberg inequality for the space $\bmo(Q_0)$. The statement of our results is given below.
  
\begin{theorem}\label{jn-blo}
	Let $Q_0$ be a finite cube in $\Rn$ and  $f\in \blo(Q_0)$ with $0<\beta\leq n$. Then for every $t>0$ and finite sub-cube $Q'\subset Q_0$, the following inequality holds
	\begin{align}\label{th3.1eq}
		\cont\left(\left\{x\in Q^\prime \, :\, f(x)-\enf_{y \in Q^\prime}f(y)>t\right\}\right)\leq c_1\exp\left(-\frac{c_2t}{\|f\|_{\blo(Q_0)}}\right)\cont(Q^\prime),
	\end{align}
	where $c_1$ and $c_2$ are positive constants independent of $Q'$ and $Q_0$.
\end{theorem}
	\begin{proof}
		Let us consider $\mathcal{D}(Q')$, the collection of all dyadic cubes generated from $Q'$, and $\contq$ is the dyadic Hausdorff content subordinate to $Q'$ (see \cite{Chen-Spector}*{Page~8} for the definition). Moreover, from  \cite{Chen-Spector}*{Proposition~2.10} we know that $\cont$ and $\contq$ are equivalent, i.e., there exists a universal constant $C>0$ (independent of $Q'$) such that for any $E\subset\Rn$ there holds
        \begin{align}\label{equi-contq}
            \frac{1}{C}\contq(E)\leq \cont(E)\leq C\contq(E).
        \end{align}
        Note that by \eqref{equi-contq}, the seminorm on $\blo(Q_0)$ is equivalent to the following seminorm defined with respect to $\contq$:
        \begin{align*}
        \|g\|_{\bloq(Q_0)}=\sup_{Q}\frac{1}{\contq(Q)}\int_{Q}\left(g(x)-\enf_{y \in Q} g(y)\right)\dcq(x),    
        \end{align*}
        where the supremum is taken over any finite subcubes $Q$ of $Q_0$ and parallel to $Q_0$.
        Therefore, to prove \eqref{th3.1eq}, it is enough to work with $\contq$ and the seminorm $\|\cdot\|_{\bloq(Q_0)}$.  Without loss of generality, we may assume $\|f\|_{\bloq(Q_0)}=1$.
        
        For any finite sub-cube $Q\subset Q_0$, we define
		\begin{align*}
			E_{Q}:=	\left\{x\in Q \, :\, f(x)-\enf_{y \in Q}f(y)>t\right\}.
		\end{align*}
		Furthermore, it follows from the monotonicity of the Choquet integral that
		\begin{align*}
			\contq(E_{Q})&\leq \frac{1}{t}\int_{E_{Q}}\left(f(x)-\enf_{y \in Q}f(y)\right)\dcq(x)\\&\leq \frac{1}{t}\int_{Q}\left(f(x)-\enf_{y \in Q}f(y)\right)\dcq(x)\\& \leq \frac{\contq(Q)}{t}\frac{1}{\contq(Q)}\int_{Q}\left(f(x)-\enf_{y \in Q}f(y)\right)\dcq(x)\\&\leq  \frac{\contq(Q)}{t} \|f\|_{\bloq(Q_0)}= \frac{\contq(Q)}{t},
		\end{align*}
holds for any cube $Q \subset Q_0$.
        
		Hence, for every $t>0$, defining $F(t)=\frac{1}{t}$, we have
		\begin{align}\label{smallest}
			\contq(E_{Q})\leq F(t) \contq ({Q}).
		\end{align}
		
Let $s>1$ and any $t\in (0,\infty)$ which satisfy $2^{\beta+1}s\leq t$. Now, for the cube $Q'\subset Q_0$, we have 
		\begin{align*}
			\frac{1}{\contq(Q')}\int_{Q'}\left(f(x)-\enf_{y \in Q'}f(y)\right)\dcq(x) \leq \|f\|_{\bloq(Q_0)}=1<s.
		\end{align*}				
		 
        Hence, applying \cite{BCRS}*{Theorem~C}, for the function $\left(f(x)-\displaystyle\enf_{y \in Q'}f(y)\right)$, there exists a countable collection of non-overlapping dyadic cubes $\{Q_k\}\subset \mathcal{D}(Q')$ with $Q_k\subset Q'$ such that:
		\begin{enumerate}
			\item $s<\frac{1}{\contq(Q_k)}\int_{Q_k}\left(f(x)-\displaystyle\enf_{y \in Q'}f(y)\right)\dcq(x) \leq 2^\beta s$ for all $k$, 
			\item $f(x)-\displaystyle\enf_{y \in Q'}f(y)\leq s$ for $\contq$-a.e. $x\in Q'\setminus \cup_k Q_k$.
		\end{enumerate} 
		An application of \cite{BCRS}*{Proposition~2.1} to this family yields a subfamily $\{ Q_{k_j}\}$ and non-overlapping ancestors $\{ \tilde{Q}_j\}$ which satisfy
		\begin{enumerate}
			\item[(a)] \begin{align*}
				\bigcup_{k} Q_k \subset \bigcup_{j} Q_{k_j} \cup \bigcup_{j} \tilde{Q}_j,
			\end{align*}
			\item[(b)] \begin{align*}
				\sum_{Q_{k_j} \subset P} \contq(Q_{k_j}) \leq 2\contq(P), \text{ for each cube } P\in \mathcal{D}(Q'),
			\end{align*}
			\item[(c)] For each $\tilde{Q}_j$, \begin{align*}
				\contq(\tilde{Q}_j) \leq \sum_{Q_{k_i} \subset \tilde{Q}_j}\contq(Q_{k_i}),
			\end{align*}
			\item[(d)] Moreover, for every nonnegative function $u\in L^1( \bigcup_{j} Q_{k_j}, \contq)$, there exists a constant $C'>0$ such that
			\begin{align*}
				\sum_{j}\int_{Q_{k_j}} u \dcq \leq C^\prime \int_{\bigcup_{j} Q_{k_j}} u\dcq .
			\end{align*}
		\end{enumerate}
		
	To avoid redundancy in this covering, if $Q_{k_j}\subset \tilde{Q}_m$ for some $j$, $m$, we will not utilize those cubes $Q_{k_j}$. The fact $s\leq 2^{-(\beta+1)} t\leq t$ and $(2)$ implies that 
		\begin{align*}
			E_{Q'}\subset  \bigcup_k Q_k \cup N \subset \bigcup_{j, Q_{k_j} \not \subset \tilde{Q}_{m}} Q_{k_j} \cup \bigcup_{j} \tilde{Q}_j \cup N,
		\end{align*}
	where $N$ is a null set with respect to the $\contq$. Now, taking the intersection with $E_{Q'}$, we will get
		\begin{align*}
			E_{Q'}&\subset \bigcup_{j, Q_{k_j} \not \subset \tilde{Q}_{m}} \left\{x\in Q_{k_j} \, :\, f(x)-\enf_{y \in Q'}f(y)>t\right\}\\& \cup \bigcup_{j} \left\{x\in \tilde{Q}_j\cap Q' \, :\, f(x)-\enf_{y \in Q'}f(y)>t\right\} \cup N,
		\end{align*}
		and so by countable subadditivity of the content $\contq$, we can estimate
		\begin{align}\label{blo-1}
			\contq(E_{Q'})&\leq \sum_{j, Q_{k_j} \not \subset \tilde{Q}_{m}}\contq\left(\left\{x\in Q_{k_j} \, :\, f(x)-\enf_{y \in Q'}f(y)>t\right\}\right)\nonumber\\&+ \sum_{j}\contq\left(\left\{x\in\tilde{Q}_j\cap Q' \, :\, f(x)-\enf_{y \in Q'}f(y)>t\right\}\right).
		\end{align}
		
Moreover, for any member of the sub-family $\{Q_{k_j}\}$, we have 
		\begin{align*}
			\int_{Q_{k_j}}\left(f(x)-\enf_{y \in Q_{k_j}}f(y)\right)\dcq(x)\leq \int_{Q_{k_j}}\left(f(x)-\enf_{y \in Q'}f(y)\right)\dcq(x).	
		\end{align*}
		Hence, by $(1)$, the sub-linearity of the Choquet integral, and the relation above, we have
		\begin{align}\label{diff-ess}
			\left|\enf_{y \in Q_{k_j}}f(y)-\enf_{y \in Q'}f(y)\right|&=\frac{1}{\contq(Q_{k_j})}\int_{Q_{k_j}}\left|\enf_{y \in Q_{k_j}}f(y)-\enf_{y \in Q'}f(y)\right|\dcq(x)\nonumber \\&\leq \frac{1}{\contq(Q_{k_j})}\int_{Q_{k_j}}\left|f(x)-\enf_{y \in Q_{k_j}}f(y)\right|\dcq(x)\nonumber \\&+\frac{1}{\contq(Q_{k_j})}\int_{Q_{k_j}}\left|f(x)-\enf_{y \in Q'}f(y)\right|\dcq(x)\nonumber\\&\leq \frac{2}{\contq(Q_{k_j})}\int_{Q_{k_j}}\left(f(x)-\enf_{y \in Q'}f(y)\right)\dcq(x)\nonumber\\& \leq 2^{\beta+1} s.
		\end{align}
		
	Therefore, using the above inequality, we estimate the first part of \eqref{blo-1} as
		\begin{align*}
			&\contq\left(\left\{x\in Q_{k_j} \, :\, f(x)-\enf_{y \in Q'}f(y)>t\right\}\right)\\& = \contq\left( \left\{x\in Q_{k_j} \, :\, f(x)-\enf_{y \in Q_{k_j}}f(y)+\enf_{y \in Q_{k_j}}f(y)-\enf_{y \in Q'}f(y)>t \right\} \right) \\& \leq \contq\left(\left\{x\in Q_{k_j} \, :\, f(x)-\enf_{y \in Q_{k_j}}f(y)>t-2^{\beta+1}s\right\}\right)\\& \leq F(t-2^{\beta+1}s)\contq(Q_{k_j}).
		\end{align*}
		
		Next, we estimate the second part of \eqref{blo-1}, which is divided into two cases. We consider the case $\tilde{Q}_j \cap Q' = Q'$. Then we have $Q' \subset \tilde{Q}_j$ and thus
		\begin{align*}
			&\contq\left(\left\{x\in \tilde{Q}_j \cap Q' \, :\, f(x)-\enf_{y \in Q'}f(y)>t\right\}\right)\\& =\contq\left(\left\{x\in Q' \, :\, f(x)-\enf_{y \in Q'}f(y)>t\right\}\right)\\& \leq \contq\left(\left\{x\in Q' \, :\, f(x)-\enf_{y \in Q'}f(y)>t-2^{\beta+1}s\right\}\right)\\&\leq  F(t-2^{\beta+1}s)\contq(Q')\\&\leq F(t-2^{\beta+1}s)\contq(\tilde{Q}_j).
		\end{align*}	
		
As we are in a dyadic setting, the only remaining case is $\tilde{Q}_j \cap Q'=\tilde{Q}_j$. Now, from the construction (c), there exists at least one $Q_{k_i}$ such that $Q_{k_i}\subset \tilde{Q}_j$. Then, using \eqref{diff-ess} and the inequality $\displaystyle\enf_{y \in Q_{k_i}}f(y)\geq \enf_{y \in \tilde{Q}_j}f(y)$, we have
		\begin{align*}
			&\contq\left(\left\{x\in \tilde{Q}_j \cap Q' \, :\, f(x)-\enf_{y \in Q'}f(y)>t\right\}\right)\\& = \contq\left( \left\{x\in \tilde{Q}_j \, :\, f(x)-\enf_{y \in Q_{k_i}}f(y)+\enf_{y \in Q_{k_i}}f(y)-\enf_{y \in Q'}f(y)>t \right\} \right) \\& \leq \contq\left(\left\{x\in \tilde{Q}_j \, :\, f(x)-\enf_{y \in Q_{k_i}}f(y)>t-2^{\beta+1}s\right\}\right)\\& \leq \contq\left(\left\{x\in \tilde{Q}_j \, :\, f(x)-\enf_{y \in \tilde{Q}_j} f(y)>t-2^{\beta+1}s\right\}\right)\\& \leq F(t-2^{\beta+1}s)\contq(\tilde{Q}_j).
		\end{align*}
		
	Finally, combining all these estimates and applying them in \eqref{blo-1}, and using (c), we obtain
		\begin{align*}
			\contq(E_{Q'})&\leq \sum_{j, Q_{k_j} \not \subset \tilde{Q}_{m}}F(t-2^{\beta+1}s)\contq(Q_{k_j})+ \sum_{j}F(t-2^{\beta+1}s)\contq(\tilde{Q}_j)\\&\leq F(t-2^{\beta+1}s) \sum_{j}\contq(Q_{k_j})\\&\leq \frac{F(t-2^{\beta+1}s)}{s} \sum_{j} \int_{Q_{k_j}}\left(f(x)-\enf_{y \in Q'}f(y)\right)\dcq(x) \\&\leq C^\prime \frac{F(t-2^{\beta+1}s)}{s}  \int_{\cup_j Q_{k_j}}\left(f(x)-\enf_{y \in Q'}f(y)\right)\dcq(x),
		\end{align*}
	where we use the fact that there is no redundancy among the cubes in the collections, and in the last two lines, we apply (1) and (d) successively. Then, using the monotonicity of the integral, the inclusion $\cup_j Q_{k_j} \subset Q'$, and the assumption $\|f\|_{\bloq(Q_0)}=1$, we have
		\begin{align*}
			\contq(E_{Q'})\leq C^\prime \frac{F(t-2^{\beta+1}s)}{s}  \int_{Q'}\left(f(x)-\enf_{y \in Q'}f(y)\right)\dcq(x)\leq  F_1(t)\contq(Q'),
		\end{align*}
		where 
		\begin{align*}
			F_1(t):= C^\prime \frac{F(t-2^{\beta+1}s)}{s}.
		\end{align*}
		
		By continuing this process inductively for any $l \geq 1$, we derive for any $2^{\beta+1}s\leq t$, there holds
		\begin{align*}
			\contq(E_{Q'})\leq  F_l(t)\contq(Q'),
		\end{align*}
		where 
		\begin{align*}
			F_l(t):= C^\prime \frac{F_{l-1}(t-2^{\beta+1}s)}{s} \quad \text{ and } \quad F_0(t)=F(t).
		\end{align*}		
		
		Now, fix $t > 0$, and suppose that
		\begin{align}\label{range}
			l\,2^{\beta+1}s<t\leq (l+1)2^{\beta+1}s,
		\end{align}
		for some $l\geq 1$, then
		\begin{align*}
			\contq(E_{Q'}) &\leq \contq\left(\left\{x\in Q' \, :\, f(x)-\enf_{y \in Q'}f(y)>l\,2^{\beta+1}s\right\}\right)\\&\leq  F_l(l\,2^{\beta+1}s)\contq(Q')\\&= C^\prime \frac{F_{l-1}(l\,2^{\beta+1}s-2^{\beta+1}s)}{s} \contq(Q')\\&={C^\prime}^{l-1}\frac{F(2^{\beta+1}s)}{s^{(l-1)}}\contq(Q')=\frac{{C^\prime}^{l-1}}{2^{\beta +1}s^{l}}\contq(Q')\\&=\frac{\contq(Q')}{2^{\beta+1}}e^{(l-1)\log C'}e^{-l\log s}.
		\end{align*}
	Now, choose $s = {C'} e > 1$, and \eqref{range} implies that
		\begin{align*}
			e^{(l-1)\log C'}e^{-l\log s}=\frac{1}{C'}e^{-l}\leq \frac{1}{C'}ee^{-\frac{t}{2^{\beta+1}{C'}e}}.
		\end{align*}
		Hence, this establishes that for any $t>2^{\beta+1}{C'}e$, the following estimate holds
		\begin{align*}
			\contq(E_{Q'})\leq C_1 e^{-C_2t}\contq(Q'),
		\end{align*}
		where 
		\begin{align*}
			C_1=\frac{e}{C'2^{\beta+1}} \quad \text{ and }\quad C_2=\frac{1}{2^{\beta+1}{C'}e}.
		\end{align*}
		
		Finally, if $t \leq 2^{\beta+1} {C'} e$, then, using the estimate,
		\begin{align*}
			\contq(E_{Q'})\leq \contq(Q')\leq e^{2^{\beta+1}{C'}e}  e^{-t} \contq(Q'),
		\end{align*}
		and choosing 
		\begin{align*}
			c_1=\max\{C_1,e^{2^{\beta+1}{C'}e}\} \quad \text{ and }\quad c_2=\min\{C_2,1\},
		\end{align*}
		the result follows for every $t>0$.
	\end{proof}

As a consequence of Theorem \ref{jn-blo}, we have the following corollary.

    \begin{corollary}
       Let $Q_0$ be a finite cube in $\Rn$ and $ 1\leq p<\infty$ and $0<\beta\leq n$. Then $\cblo^{\beta,p}(Q_0)=\blo(Q_0)$, with respect to the corresponding norms. 
    \end{corollary}
	\begin{proof}
		Let $1<p<\infty$. Assume $f \in \blo(Q_0)$ and take $Q' \subset Q_0$. Apply Theorem~\ref{jn-blo} to the cube $Q'$, along with a change of variables and Hölder's inequality, we get
		\begin{align*}
			&\frac{1}{\cont(Q')}\int_{Q^\prime}\left(f(x)-\enf_{y \in Q^\prime} f(y)\right)^p\dc(x)\\
            &= \frac{1}{\cont(Q')}\int_0^\infty \cont \left(\{x\in Q'\,:\, \left(f(x)-\enf_{y \in Q^\prime} f(y)\right)^p>t\}\right)\dt\\
            &= \frac{p}{\cont(Q')}\int_0^\infty \lambda^{p-1}\cont \left(\{x\in Q'\,:\, \left(f(x)-\enf_{y \in Q^\prime} f(y)\right)>\lambda\}\right){\rm d}\lambda\\
            &\leq p c_1 \int_0^\infty \lambda^{p-1} e^{-c_2 \lambda/\|f\|_{\blo(Q_0)}} \,{\rm d}\lambda= C\|f\|_{\blo(Q_0)}^p,
		\end{align*}
		and by considering the supremum, we conclude that
		\begin{align*}
			\blo(Q_0) \subset\cblo^{\beta,p}(Q_0).
		\end{align*}
	The reverse inclusion follows by applying Hölder's inequality with exponents $\frac{1}{p} + \frac{1}{p/(p-1)} = 1$. This completes the proof.
	\end{proof} 
As an application of Theorem~\ref{jn-blo}, we obtain the following corollary, which establishes the inclusion relationships among the $\blo$ spaces for various values of $\beta$ with $0 < \beta \leq n$.  
\begin{corollary}\label{nesting-blo}
	Let $0 < \gamma \leq \beta \leq n$, and suppose that $f \in \blob{\gamma}(Q_0)$. Then $f \in \blo(Q_0)$, and
	\begin{align*}
		\|f\|_{\blo(Q_0)} \leq C\|f\|_{\blob{\gamma}(Q_0)},
	\end{align*}
	for a positive constant $C = C(\gamma, \beta, n)$ independent of $u$.
\end{corollary}
\begin{proof}
	Let $Q' \subset Q_0$. Then, by applying \cite{ChenOoiSpector}*{Lemma~2.2}, Theorem~\ref{jn-blo}, and \eqref{equiv-cont-cube}, we obtain
    \begin{align*}
			&\frac{1}{\cont(Q')}\int_{Q^\prime}\left(f(x)-\enf_{y \in Q^\prime} f(y)\right)\dc(x)\\
            &= \frac{1}{\cont(Q')}\int_0^\infty \cont \left(\{x\in Q'\,:\, \left(f(x)-\enf_{y \in Q^\prime} f(y)\right)>t\}\right)\dt\\
            &\leq  \frac{C}{\cont(Q')}\int_0^\infty \left[\contb{\gamma} \left(\{x\in Q'\,:\, \left(f(x)-\enf_{y \in Q^\prime} f(y)\right)>t\}\right)\right]^{\frac{\beta}{\gamma}}\dt \\
            &\leq  \frac{C}{\cont(Q')}\int_0^\infty \left[c_1\exp\left(-\frac{c_2t}{\|f\|_{\blob{\gamma}(Q_0)}}\right)\contb{\gamma}(Q^\prime)\right]^{\frac{\beta}{\gamma}}\dt\\
            &=C\frac{\left(\contb{\gamma}(Q')\right)^{\frac{\beta}{\gamma}}}{\cont(Q')}c_1^{\frac{\beta}{\gamma}}\int_0^\infty \left[\exp\left(-\frac{c_2\beta t}{\gamma\|f\|_{\blob{\gamma}(Q_0)}}\right)\right]^{\frac{\beta}{\gamma}}\dt\leq C\|f\|_{\blob{\gamma}(Q_0)},
		\end{align*}
        and by taking the supremum over all $Q' \subset Q_0$, the result follows.
\end{proof}
	
\section{Bounded Mean oscillation spaces and fractional maximal function \texorpdfstring{$\fma$}{fma}}\label{Sec:BMO-improve}
In this section, we will prove our first main result, Theorem \ref{MainThm2}. To prove this, we need the following known result about the boundedness of $\fma$ from $L^p(\Rn)$ to $L^p(\cont, \Rn)$ for $\beta=n-\alpha p$.
\begin{lemma}[\cite{Sawyer}]\label{Adamslaylemma} 
Let $0 < \alpha < n.$ If $1 < p < \frac{n}{\alpha}$ and $\beta = n - \alpha p$, then there exists a constant $C$ depending on $n$, $\alpha$, and $p$ such that \begin{align*} \left( \int_{\Rn} (\fma f)^p \dcb{n- \alpha p} \right)^{\frac{1}{p}} \leq C \left( \int_{\Rn} |f|^p  \dx \right)^{\frac{1}{p}}. 
\end{align*} 
\end{lemma}

The proof of Lemma~\ref{Adamslaylemma} is based on the one-weight inequality established in \cite{Sawyer}*{Theorem B}. Furthermore, it is presented in \cite{AdamsChoquet}*{Theorem 7}, where the result is discussed in the framework of Hausdorff content, as stated in Lemma~\ref{Adamslaylemma}.

Now, we give the proof of Theorem \ref{MainThm2}.
\begin{proof}[Proof of Theorem \ref{MainThm2}]
By Corollary~\ref{nesting-blo}, it suffices to establish the theorem in the case $0 < \beta < n - \alpha$. Furthermore, by Lemma~\ref{abs-bmo} and the definition of $\mathcal{M}_{\alpha}$, we may, without loss of generality, assume that $f$ is nonnegative. If $f$ is not nonnegative, we may replace it with $|f|$, since the fractional maximal function $\mathcal{M}_{\alpha}$ remains unaffected by taking absolute values.

Let $Q$ be a cube such that $Q \subset Q_0$. Choose $p$ with $1 < p < \frac{n}{\alpha}$ such that $\beta = n - \alpha p$. The existence of such a $p$ is guaranteed by the assumption $\beta < n - \alpha$. Applying Lemma~\ref{diffBMO} with $\gamma = 0$, we obtain a cube $P_{2Q}$ satisfying 
\begin{enumerate}[label=(\roman*)]
    \item $P_{2Q} \subseteq Q_0 \cap 4Q$,
    \item $P_{2Q} \supseteq Q_0 \cap 2Q$,
    \item The side length of $P_{2Q}$ satisfies
    \[
    \ell(P_{2Q}) = \begin{cases}
        \ell(2Q) & \text{if } \ell(2Q) < \ell(Q_0), \\[6pt]
        \ell(Q_0) & \text{if } \ell(2Q) \geq \ell(Q_0),
    \end{cases}
    \]
    \item The following inequality holds:
    \begin{align*}
    \left(\frac{1}{\ell(2Q)^n} \int_{2Q \cap Q_0} |f(y) - f_{P_{2Q}}|^p \, \dy \right)^{\frac{1}{p}}
    \leq \left( \frac{1}{\ell(P_{2Q})^n} \int_{P_{2Q}} |f(y) - f_{P_{2Q}}|^p \, \dy\right)^{\frac{1}{p}}.
    \end{align*}
\end{enumerate}
We now decompose $f$ as $f = g + h$, where
\begin{align*}
   g := (f - f_{P_{2Q}})\chi_{(2Q \cap Q_0)} \quad \text{and} \quad h := f - g = f_{P_{2Q}}\chi_{(2Q \cap Q_0)} + f\chi_{(Q_0 \setminus 2Q)}. 
\end{align*}
Applying H\"older's inequality, Lemma~\ref{Adamslaylemma}, and \eqref{equiv-cont-cube}, we deduce that there exists a constant $C = C(n,\alpha,\beta) > 0$ such that
\begin{align}\label{estimateg}
    \frac{1}{\cont(Q)} \int_Q \mathcal{M}_\alpha g(y) \dc(y)
    &\leq \frac{1}{\cont(Q)^{1/p}} \left( \int_{\mathbb{R}^n} (\mathcal{M}_\alpha g(y))^p \, \dc(y) \right)^{\!1/p} \nonumber\\[6pt]
    &\leq \frac{C}{\cont(Q)^{1/p}} \left( \int_{\mathbb{R}^n} |g(y)|^p  \dy \right)^{\!1/p} \nonumber\\[6pt]
    &= \frac{C}{\cont(Q)^{1/p}} \left( \int_{2Q \cap Q_0} |f(y) - f_{P_{2Q}}|^p \, \dy \right)^{\!1/p} \nonumber\\[6pt]
\nonumber    &\leq C \frac{\ell(2Q)^{n/p}}{\ell(Q)^{\beta/p}} \left( \frac{1}{\ell(2Q)^n} \int_{2Q \cap Q_0} |f(y) - f_{P_{2Q}}|^p \dy \right)^{\!1/p}\\
  &\leq  C 2^{\frac{n}{p}}  \ell(Q_0)^\alpha \left( \frac{1}{\ell(2Q)^n} \int_{2Q \cap Q_0} |f(y) - f_{P_{2Q}}|^p \dy \right)^{\!1/p},
\end{align}
where in the last inequality, we have used the fact $n/p - \beta/p = \alpha$ and $Q \subset Q_0$.
By property (iii) of the cube $P_{2Q}$, we have
\[
\left( \frac{1}{\ell(2Q)^n} \int_{2Q \cap Q_0} |f(y)-f_{P_{2Q}}|^p \dy \right)^{\!1/p} \leq  \left( \frac{1}{\ell(P_{2Q})^n} \int_{P_{2Q}} |f(y)-f_{P_{2Q}}|^p \dy \right)^{\!1/p}.
\]
 Consequently,
\begin{align*}
    \frac{1}{\cont(Q)} \int_Q \mathcal{M}_\alpha g(y) \dc(y)
    &\leq C 2^{\frac{n}{p}}  \ell(Q_0)^\alpha \left( \frac{1}{\ell(P_{2Q})^n} \int_{P_{2Q}} |f(y)-f_{P_{2Q}}|^p \dy \right)^{\!1/p} \\[6pt]
    &\leq C 2^{\frac{n}{p}}  \ell(Q_0)^\alpha \|f\|_{\cbmo(Q_0)},
\end{align*}
where the last step follows from the property (i) of $P_{2Q}$ together with the definition of the \(\cbmo(Q_0)\)-norm.

Next, we establish an upper bound for $\mathcal{M}_\alpha h(x_0)$ for each $x_0 \in Q$. Fix $x_0 \in Q$ and consider a cube $\tilde{Q} \subseteq Q_0$ with $x_0 \in \tilde{Q}$.

\textbf{Case:} $\tilde{Q} \subseteq 2Q$.

If $\tilde{Q}$ is entirely contained in $2Q$, then $\ell(\tilde{Q}) \leq 2\ell(Q)$. Employing properties (ii) and (iii) of $P_{2Q}$, together with the definition of the fractional maximal function $\mathcal{M}_\alpha$, we deduce that the estimate
\begin{align}\label{hassumption1}
	\frac{1}{\ell(\tilde{Q})^{n-\alpha}} \int_{\tilde{Q}} |h| \dy 
	&= \ell(\tilde{Q})^\alpha f_{P_{2Q}} \nonumber \\
    \nonumber    & =  \left\{
\begin{array}{ll}
\ell(\tilde{Q})^\alpha \left(  \frac{1}{\ell(2Q)^n }  \int_{P_{2Q}}  f(y) \dy \right) & \text{, if } \ell(2Q) < \ell(Q_0) \\
\ell(\tilde{Q})^\alpha \left( \frac{1}{\ell(Q_0)^n }  \int_{P_{2Q}}  f (y)\dy \right)  & \text{, if } \ell(2Q) \geq \ell(Q_0) 
\end{array}
\right. \\
\nonumber & \leq    \left\{
\begin{array}{ll}
\frac{1}{\ell(2Q)^{n-\alpha }} \int_{P_{2Q}}  f (y)\dy  & \text{, if } \ell(2Q) < \ell(Q_0) \\
\frac{1}{\ell(Q_0)^{n-\alpha}} \int_{P_{2Q}}  f (y)\dy  & \text{, if } \ell(2Q) \geq \ell(Q_0) 
\end{array}
\right. \\
\nonumber& =   \frac{1}{\ell(P_{2Q})^{n-\alpha}} \int_{P_{2Q}}  f (y)\dy   \\
	&\leq \mathcal{M}_\alpha f(z)
\end{align}
holds for every $z \in Q$.

\textbf{Case:} $\tilde{Q} \cap (2Q)^c \neq \emptyset$.

If the cube $\tilde{Q}$ is not fully contained in $2Q$, then the cube $16\tilde{Q}$ contains both $\tilde{Q}$ and $4Q$. By applying Lemma~\ref{diffBMO}, there exists a cube $P_{16\tilde{Q}}$ such that
\begin{enumerate}[label=(\roman*')]
    \item $P_{16\tilde{Q}} \subseteq Q_0 \cap 32\tilde{Q}$,
    \item $P_{16\tilde{Q}} \supseteq Q_0 \cap 16\tilde{Q}$,
    \item The side length of $P_{16\tilde{Q}}$ satisfies
    \begin{align*}
       \ell(P_{16\tilde{Q}})=\begin{cases}
    \ell(16\tilde{Q}) & \text{if } \ell(16\tilde{Q}) < \ell(Q_0), \\[6pt]
    \ell(Q_0) & \text{if } \ell(16\tilde{Q}) \geq \ell(Q_0),
    \end{cases} 
    \end{align*}
    \item The following inequality holds:
    \begin{align*}
          \frac{1}{\ell(16\tilde{Q})^n} \int_{16\tilde{Q} \cap Q_0} |f(y) - f_{P_{16\tilde{Q}}}| \dy 
    \leq  \frac{1}{\ell(P_{16\tilde{Q}})^n} \int_{P_{16\tilde{Q}}} |f(y) - f_{P_{16\tilde{Q}}}| \dy.
    \end{align*}
\end{enumerate}
We now split the integral as follows
\begin{align*}
    \frac{1}{\ell(\tilde{Q})^{n-\alpha}}\int_{\tilde{Q}} |h(y)| \dy 
    &\leq \ell(\tilde{Q})^\alpha \left(\frac{1}{\ell(\tilde{Q})^n} \int_{\tilde{Q}} |h(y) - f_{P_{16\tilde{Q}}}| \dy \right) \\
    &\quad + \ell(\tilde{Q})^\alpha \left(\frac{1}{\ell(P_{16\tilde{Q}})^n} \int_{P_{16\tilde{Q}}} f(y) \dy \right) \\
    &\leq  2^{4n} \ell(\tilde{Q})^\alpha \left(\frac{1}{\ell(16\tilde{Q})^n} \int_{16\tilde{Q} \cap Q_0} |h(y) - f_{P_{16\tilde{Q}}}| \dy \right) \\
    &\quad + \ell(\tilde{Q})^\alpha \left(\frac{1}{\ell(P_{16\tilde{Q}})^n} \int_{P_{16\tilde{Q}}} f(y) \dy \right) \\
    &:= I + II.
\end{align*}

For $I$, we note that the definition $h = f_{P_{2Q}} \chi_{(2Q \cap Q_0)} + f\chi_{(Q_0 \setminus 2Q)}$ implies
\begin{align}\label{sqsq}
 \nonumber    \frac{1}{ \ell(16 \tilde{Q})^n} \int_{16 \tilde{Q} \cap Q_0} |h(y) - f_{P_{16\tilde{Q}}}| \dy
 &= \frac{1}{ \ell(16 \tilde{Q})^n} \int_{2Q \cap Q_0} |h(y) - f_{P_{16\tilde{Q}}}| \dy \\
 \nonumber    &  + \frac{1}{ \ell(16 \tilde{Q})^n} \int_{(16 \tilde{Q} \setminus 2Q) \cap Q_0} |h(y) - f_{P_{16\tilde{Q}}}| \dy \\
 \nonumber     &= \frac{1}{ \ell(16 \tilde{Q})^n} \int_{2Q \cap Q_0} |f_{P_{2Q}} - f_{P_{16\tilde{Q}}}| \dy \\
     & + \frac{1}{ \ell(16 \tilde{Q})^n} \int_{(16 \tilde{Q} \setminus 2Q) \cap Q_0} |f(y) - f_{P_{16\tilde{Q}}}| \dy.
\end{align}
 We now estimate the first term from the last line of (\ref{sqsq}). By property (i) of $P_{2Q}$ and property (ii') of $P_{16 \tilde{Q}}$, we have
$$P_{16 \tilde{Q}} \supseteq Q_0 \cap 16 \tilde{Q} \supseteq Q_0 \cap 4Q \supseteq P_{2Q}.$$
From this chain of inclusions, it follows that
\begin{align*}
    | f_{P_{2Q}} - f_{P_{16 \tilde{Q}}} |  \leq \frac{1}{ \ell(P_{2Q})^n} \int_{P_{2Q}} | f(y) - f_{P_{16 \tilde{Q}}}| \dy 
     \leq \frac{1}{\ell(P_{2Q})^n} \int_{P_{16 \tilde{Q}}} |f(y) - f_{P_{16 \tilde{Q}}}| \dy.
\end{align*}
Applying this estimate, we obtain
\begin{align}\label{squar1}
   \frac{1}{\ell(16 \tilde{Q})^n } \int_{2Q \cap Q_0} |f_{P_{2Q}} - f_{P_{16} \tilde{Q}}| \dy 
     \leq \frac{|2Q \cap Q_0|}{ \ell(P_{2Q})^n} \frac{1}{\ell(16 \tilde{Q})^n} \int_{P_{16 \tilde{Q}}} |f(y) - f_{P_{16 \tilde{Q}}}| \dy 
     \leq \Vert f \Vert_{BMO(Q_0)},
\end{align}
where in the last step we have used the inequality $\ell(16 \tilde{Q}) \geq \ell(P_{16 \tilde{Q}})$ (see (iii')) and the estimate
\begin{align*}
 \frac{|2Q \cap Q_0|}{\ell(P_{2Q})^n} = \begin{cases}
\frac{|2Q \cap Q_0|}{\ell(2Q)^n} & \text{if } \ell(2Q) < \ell(Q_0) \\[6pt]
\frac{|2Q \cap Q_0|}{\ell(Q_0)^n} & \text{if } \ell(2Q) \geq \ell(Q_0)
\end{cases}   
\end{align*}
which in either case is at most $1$. Next, we turn to the second term from the last line of \eqref{sqsq}. By property (iv') of $P_{16 \tilde{Q}}$, we have
\begin{align}\label{squar2}
 \nonumber   \frac{1}{ \ell(16 \tilde{Q})^n} \int_{(16 \tilde{Q} \setminus 2Q) \cap Q_0} |f(y) - f_{P_{16\tilde{Q}}}| \dy 
\nonumber   & \leq  \frac{1}{ \ell(16 \tilde{Q})^n} \int_{16 \tilde{Q} \cap Q_0} |f(y) - f_{P_{16\tilde{Q}}}| \dy \\
\nonumber   & \leq  \frac{1 }{ \ell(P_{16 \tilde{Q}})^n} \int_{P_{16 \tilde{Q}}}  | f(y) -f_{P_{16 \tilde{Q}}}| \dy \\
   & \leq  \Vert f \Vert_{\cbmo(Q_0)}.
\end{align}
Combining (\ref{sqsq}), (\ref{squar1}), and (\ref{squar2}), we conclude that
\begin{align}\label{st1}
   \nonumber I &\leq 2^{4 n +1}  \ell(\tilde{Q})^\alpha \Vert f \Vert_{\cbmo(Q_0)}\\
    & \leq 2^{4 n +1} \ell(Q_0)^\alpha \Vert f \Vert_{\cbmo(Q_0)}.
\end{align}

For $II$, we observe that by the inclusion $\tilde{Q} \subseteq Q_0$ and property (iii') of $P_{16 \tilde{Q}}$, we have $\ell(\tilde{Q}) \leq \ell(P_{16 \tilde{Q}})$. Thus, for every $z \in Q$,
\begin{align}\label{st2}
   II= \ell(\tilde{Q})^\alpha \left( \frac{1}{\ell(P_{16\tilde{Q}})^n} \int_{P_{16 \tilde{Q}}} f(y) \dy \right) 
  \leq \frac{1}{\ell(P_{16 \tilde{Q}})^{n-\alpha}} \int_{P_{16 \tilde{Q}}} f(y) \dy \leq \mathcal{M}_\alpha f(z),
\end{align}
where in the last inequality we have used the inclusions $$P_{16\tilde{Q}} \supseteq 16 \tilde{Q} \cap Q_0 \supseteq 4Q \cap Q_0 \supseteq Q.$$ The combination of \eqref{st1} and \eqref{st2} then gives
\begin{align}\label{hassumption2}
	\frac{1}{\ell(\tilde{Q})^{n-\alpha}} \int_{\tilde{Q}} |h| \dy 
    &\leq 2^{4n+1} \ell(Q_0)^\alpha \Vert f \Vert_{\cbmo(Q_0)} + \mathcal{M}_\alpha f(z)
\end{align}
for every $z \in Q$. 

Taking the supremum over all such $\tilde{Q}$ in the inequalities \eqref{hassumption1} and \eqref{hassumption2}, and then applying part \eqref{ywl-mf} of Remark~\ref{equi-bmo}, we conclude that for every cube $Q \subseteq Q_0$,
\begin{align*}
    \mathcal{M}_\alpha h (x_0) 
    \leq 2^{4n+1} \ell(Q_0)^\alpha \Vert f \Vert_{\cbmo(Q_0)} + \mathcal{M}_\alpha f(z)
\end{align*}
holds for every $x_0, z \in Q$. Consequently, for every cube $Q \subseteq Q_0$,
\begin{align}\label{hassumption4}
    \mathcal{M}_\alpha h (x_0) 
    \leq 2^{4n+1} \ell(Q_0)^\alpha \Vert f \Vert_{\cbmo(Q_0)} 
    + \enf_{z \in Q } \mathcal{M}_\alpha f(z)
\end{align}
for every $x_0 \in Q$. By the sublinearity of the fractional maximal function $\mathcal{M}_\alpha$, together with \eqref{estimateg} and \eqref{hassumption4}, we obtain
\begin{align*}
   & \frac{1}{\cont(Q)} \int_{Q \cap Q_0}  \left(\mathcal{M}_\alpha f (y)- \enf_{z \in Q } \mathcal{M}_\alpha f (z) \right)\dc(y)\\
    \leq& \frac{1}{\cont(Q)} \int_{Q \cap Q_0} \left(\mathcal{M}_\alpha g(y) 
    + \mathcal{M}_\alpha h (y ) - \enf_{z \in Q } \mathcal{M}_\alpha f (z) \right)\dc(y) \\
   \leq & \frac{1}{\cont(Q)} \int_{Q \cap Q_0} \left(\mathcal{M}_\alpha g(y) 
    +  \left( 2^{4n+1} \ell(Q_0)^\alpha \Vert f \Vert_{\cbmo(Q_0)} \right) \right)\dc(y) \\
    \leq &\frac{1}{\cont(Q)} \int_{Q} \mathcal{M}_\alpha g(y) \dc(y) 
    + 2^{4n+2} \ell(Q_0)^\alpha \Vert f \Vert_{\cbmo(Q_0)} \\
    \leq &\left( C 2^{\frac{n}{p}+1} + 2^{4n+2} \right) \ell(Q_0)^\alpha  \Vert f \Vert_{\cbmo (Q_0)}.
\end{align*}
Finally, using \eqref{ywl-blo}, we get the desired estimate.
\end{proof}

\section{\texorpdfstring{$\beta$}{Beta}-dimensional Bounded mean oscillation spaces and maximal function \texorpdfstring{$\mathcal{M}^{\beta}$}{} }\label{Sec:beta-BMO-Maximal function}
In this section, we prove our second main result, Theorem~\ref{asthe}. Actually, we prove more general results concerning the boundedness of the $\beta$-dimensional uncentered maximal functions on BMO spaces with respect to the dyadic Hausdorff content $\cont$, which is stated in Theorem \ref{MainThm}. While the proof follows a similar strategy to that of Theorem~\ref{MainThm2}, a key difference is that, instead of applying Lemma~\ref{Adamslaylemma}, we utilize the following $L^p$-boundedness result for $\mathcal{M}^{\beta}$, which is a consequence of \cite{ChenOoiSpector}*{Theorem~1.2, Corollary~1.4}. We need the following two lemmas to prove our main results. 

\begin{lemma}\label{modf-lem}
  Suppose $0<\gamma\leq \beta\leq n$, and $p\in\big(\frac{\gamma}{\beta}, \infty\big)$. Then there exists a constant $C>0$, depending only on  $\gamma$, $\beta$, $n$, and $p$ such that
  \begin{align*}
      \int_{\Rn}\left(\mathcal{M}^\beta f\right)^{p}\dcb{\gamma}\leq C \int_{\Rn}|f|^p\dcb{\gamma},
  \end{align*}
  for all $f\in L^p(\contb{\gamma}, \Rn)$.
\end{lemma}
\begin{proof}
    The proof follows from \cite{ChenOoiSpector}*{Theorem~1.2} for $\gamma=\beta$ and from \cite{ChenOoiSpector}*{Corollary~1.4} for $\gamma<\beta$ and using \eqref{equiv-cont}. 
 \end{proof}

 \begin{lemma}\label{compare-beta}
     Let $0<\beta_1\leq\beta_2\leq n$. Then, for any nonnegative function $f\in \bmob{\beta_1}(\Rn)$ and for any finite cube $Q$, the following estimate holds \begin{align*} 
     \left|\avgdb{Q}{\beta_1}-\avgdb{Q}{\beta_2}\right|\leq C \|f\|_{\bmob{\beta_1}(\Rn)}, 
     \end{align*} 
     where $C$ is a positive constant depending only on $n$, $\beta_1$, and $\beta_2$.
 \end{lemma}
 \begin{proof}
  The case $\beta_1=\beta_2$ follows trivially. Therefore, we assume $\beta_1<\beta_2$. Applying Lemma~\ref{tri_cont} and \cite{ChenOoiSpector}*{Lemma~2.2}, we have
  \begin{align*}
         \left|\avgdb{Q}{\beta_1}-\avgdb{Q}{\beta_2}\right|&\leq \frac{1}{\contb{\beta_2}(Q)}\int_{Q}|f-\avgdb{Q}{\beta_1}|\dcb{\beta_2}\\&\leq \frac{C}{\contb{\beta_2}(Q)}\left[\int_{Q}|f-\avgdb{Q}{\beta_1}|^\frac{\beta_1}{\beta_2}\dcb{\beta_1}\right]^{\frac{\beta_2}{\beta_1}}\\&\leq \frac{C}{\contb{\beta_2}(Q)}\left(\contb{\beta_1}(Q)\right)^{\frac{\beta_2-\beta_1}{\beta_1}}\left(\int_{Q}|f-\avgdb{Q}{\beta_1}|\dcb{\beta_1}\right)\\&\leq \frac{C}{\contb{\beta_2}(Q)}\left(\contb{\beta_1}(Q)\right)^{\frac{\beta_2}{\beta_1}}\left(\frac{1}{\contb{\beta_1}(Q)}\int_{Q}|f-\avgdb{Q}{\beta_1}|\dcb{\beta_1}\right)\\&\leq C\|f\|_{\bmob{\beta_1}(\Rn)}.
     \end{align*}
To conclude the proof of the lemma, in the middle, we applied H\"older's inequality with $\frac{1}{\beta_2/\beta_1}+\frac{1}{\beta_2/(\beta_2-\beta_1)}=1$ and \eqref{equiv-cont-cube} in the end.
 \end{proof}

Now, our main result of this section can be stated as follows:
\begin{theorem}\label{MainThm}
Let $0 < \beta_1 \leq  \beta_2 \leq n$.  If $f\in \bmob{\beta_1}(\Rn)$ and $\mathcal{M}^{\beta_2} f$ is not identically infinite, then there exists a constant $C > 0$, depending only on $\beta_1$, $\beta_2$, and $n$ such that
\begin{align*}
    \|\mathcal{M}^{\beta_2} f\|_{\blob{\beta_1}(\Rn)}\leq C \|f\|_{\bmob{\beta_1}(\Rn)}.
\end{align*}
\end{theorem} 
 \begin{proof}
 Let us assume $f \in \bmob{\beta_1}(\Rn)$. By the Lemma \ref{abs-bmo}, we conclude that $|f|$  also in  $\bmob{\beta_1}(\Rn)$ and $\left\| |f|\right\|_{\bmob{\beta_1}(\Rn)}\leq C \|f\|_{\bmob{\beta_1}(\Rn)}$. Therefore, since $\mathcal{M}^{\beta_2}(f) = \mathcal{M}^{\beta_2}(|f|)$, to prove the Theorem~\ref{MainThm}, it is enough to assume that $f$ is a nonnegative function. 
 
Now fix some $p\in(\frac{\beta_1}{\beta_2},\infty)$ and $p> 1$. We also fix a cube $Q' \subset \Rn$. We decompose $f= g+h$, where $g= (f-\avgdb{2Q'}{\beta_2} ) \chi_{2Q'}$ and $h=f-g= (\avgdb{2Q'}{\beta_2})\chi_{2Q'} + f \chi_{\Rn\setminus 2Q'}$. 

First, we estimate for  $\mathcal{M}^{\beta_2}g$. Using first H\"older's inequality with exponent $p>1$, and Lemma~\ref{modf-lem}, we deduce
\begin{align}\label{best-for-g-frac}
     \nonumber \frac{1}{\contb{\beta_1}(Q')} \int_{Q'} \mathcal{M}^{\beta_2}g \dcb{\beta_1} & \leq \contb{\beta_1}(Q')^{-\frac{1}{p}}\left(\int_{Q'}\left(\mathcal{M}^{\beta_2} g\right)^{p}\dcb{\beta_1}\right)^{\frac{1}{p}}\\ & \nonumber\leq  C \contb{\beta_1}(Q')^{-\frac{1}{p}}\left(\int_{\Rn}|g|^p\dcb{\beta_1}\right)^{\frac{1}{p}} \\\nonumber & \leq  C \left( \frac{1}{\contb{\beta_1}(2Q')} \int_{2Q'} |f-\avgdb{2Q'}{\beta_2}|^p \dcb{\beta_1} \right)^{1/p}\\\nonumber & \leq  C \left( \frac{1}{\contb{\beta_1}(2Q')} \int_{2Q'} |f-\avgdb{2Q'}{\beta_1}|^p \dcb{\beta_1} \right)^{1/p} + C \left|\avgdb{2Q'}{\beta_1}-\avgdb{2Q'}{\beta_2}\right| \\ & \leq C \|f\|_{\bmob{\beta_1}(\Rn)}.
		\end{align}
In the third inequality, we used the doubling property of $\contb{\beta_1}$, and in the last inequality, we used Lemma~\ref{compare-beta}.

Next, we estimate $\mathcal{M}^{\beta_2}(h)(z)$ for each $z\in Q'$. Fix $z \in Q'$ and choose another cube $\bar{Q}$ such that $z\in \bar{Q}$. If $\bar{Q} \subseteq 2Q'$, then we have $\ell(\bar{Q}) \leq 2 \ell(Q')$. In this case, we get
		\begin{align}\label{esti-h-frac-1}
			\frac{1}{\contb{\beta_2}(\bar{Q})} \int_{\bar{Q}} |h| \dcb{\beta_2}  =   \avgdb{2Q'}{\beta_2}
              \leq    \mathcal{M}^{\beta_2}(f)(x)
		\end{align}
		for any $x \in Q'$.

        Now consider the remaining case $\bar{Q}\cap (2Q')^c\neq \emptyset$. In this case $16\bar{Q}$ contains $4Q'\cup \bar{Q}$. Using the sublinearity of the integral and doubling the property of the dyadic content, we have
\begin{align*}
    \frac{1}{\contb{\beta_2}(\bar{Q})}\int_{\bar{Q}} |h| \dcb{\beta_2} 
    &\leq \frac{1}{\contb{\beta_2}(\bar{Q})} \int_{\bar{Q}} |h - \avgdb{16\bar{Q}}{\beta_2}| \dcb{\beta_2}  +  \frac{1}{\contb{\beta_2}(16\bar{Q})} \int_{16\bar{Q}} f \dcb{\beta_2} \\
    &\leq C \left(\frac{1}{\contb{\beta_2}(16\bar{Q})} \int_{16\bar{Q}} |h - \avgdb{16\bar{Q}}{\beta_2}| \dcb{\beta_2} \right) +\frac{1}{\contb{\beta_2}(16\bar{Q})} \int_{16\bar{Q}} f \dcb{\beta_2}\\
    &:= I + II.
\end{align*}

For $I$, we note that the definition $h =  (\avgdb{2Q'}{\beta_2} ) \chi_{2Q'} + f \chi_{\Rn\setminus 2Q'}$ implies
\begin{align}\label{bsqsq_bq}
    \int_{16\bar{Q}} |h - \avgdb{16\bar{Q}}{\beta_2}| \dcb{\beta_2} 
= \int_{2Q' } |\avgdb{2Q'}{\beta_2} - \avgdb{16\bar{Q}}{\beta_2}| \dcb{\beta_2} 
       +   \int_{(16 \bar{Q} \setminus 2Q') } |f- \avgdb{16\bar{Q}}{\beta_2}| \dcb{\beta_2} .
\end{align}

We now estimate the first term from the last line of (\ref{bsqsq_bq}). Using $2Q'\subset 4Q'\subset 16 \bar{Q}$, Lemma~\ref{tri_cont}, and Theorem~\ref{bmo-remark-def}, we deduce
\begin{align}\label{bsquar1}
    \frac{1}{\contb{\beta_2}(16\bar{Q})} \int_{2Q'} |\avgdb{2Q'}{\beta_2} - \avgdb{16\bar{Q}}{\beta_2}| \dcb{\beta_2} \nonumber&=  \frac{\contb{\beta_2}(2Q')}{\contb{\beta_2}(16\bar{Q})}|\avgdb{2Q'}{\beta_2} - \avgdb{16\bar{Q}}{\beta_2}|\nonumber\\& \leq \frac{1}{\contb{\beta_2}(16\bar{Q})}\int_{2Q'}|f-\avgdb{16\bar{Q}}{\beta_2}|\dcb{\beta_2}\nonumber\\
    &\leq \frac{1}{\contb{\beta_2}(16\bar{Q})}\int_{16\bar{Q}}|f-\avgdb{16\bar{Q}}{\beta_2}|\dcb{\beta_2} \nonumber \\& \leq C  \|f\|_{\bmob{\beta_2}(\Rn)}.
\end{align}

Next, we turn to the second term from the last line of \eqref{bsqsq_bq}. By using the containment $(16 \bar{Q} \setminus 2Q')\subset 16\bar{Q}$, and Theorem~\ref{bmo-remark-def}, we have
\begin{align}\label{bsquar2}
    \frac{1}{\contb{\beta_2}(16\bar{Q})}  \int_{(16 \bar{Q} \setminus 2Q') } |f- \avgdb{16\bar{Q}}{\beta_2}| \dcb{\beta_2}   \leq C \|f\|_{\bmob{\beta_2}(\Rn)}
\end{align}

Now, combining \eqref{bsquar1}, and \eqref{bsquar2} we conclude that
\begin{align}\label{bst1}
    I \leq C \|f\|_{\bmob{\beta_2}(\Rn)}.
\end{align}
For $II$, we observe that $z\in Q'\subset 16\bar{Q}$, we obtain
\begin{align}\label{bst2}
  II\leq  \mathcal{M}^{\beta_2}f(z).
\end{align}
Hence, the combination of \eqref{bst1} and \eqref{bst2} gives
\begin{align}\label{bhassumption2}
	 \frac{1}{\contb{\beta_2}(\bar{Q})}\int_{\bar{Q}}|h|\dcb{\beta_2}
    &\leq C \|f\|_{\bmob{\beta_2}(\Rn)} + \mathcal{M}^{\beta_2}f(z),
\end{align}
for every $z \in Q'$ and $C$ is some universal constant independent on $\bar{Q}$ and $Q'$. Taking the supremum over all such $\bar{Q}$ in the inequalities \eqref{bhassumption2}, we conclude that for every cube $Q'$,
\begin{align*}
    \mathcal{M}^{\beta_2} h (x_0) 
    \leq C \|f\|_{\bmob{\beta_2}(\Rn)} + \mathcal{M}^{\beta_2} f(z)
\end{align*}
holds for every $x_0, z \in Q'$. Consequently, for every cube $Q'$, we have
\begin{align}\label{bhassumption4}
     \mathcal{M}^{\beta_2} h (x_0) 
    \leq C \|f\|_{\bmob{\beta_2}(\Rn)}
    + \esnf_{z \in Q' } \mathcal{M}^{\beta_2} f (z)
\end{align}
for every $x_0 \in Q'$. By the sublinearity of the fractional maximal function $\mathcal{M}^{\beta_2}$, together with \eqref{best-for-g-frac} and \eqref{bhassumption4}, we obtain
\begin{align*}
   & \frac{1}{\contb{\beta_1}(Q')} \int_{Q'}  \left(\mathcal{M}^{\beta_2} f (y)- \esnf_{z \in Q' }\mathcal{M}^{\beta_2} f (z) \right)\dcb{\beta_1}(y)\\
    \leq& \frac{1}{\contb{\beta_1}(Q')}  \int_{Q' } \left(\mathcal{M}^{\beta_2} g(y) 
    + \mathcal{M}^{\beta_2} h (y ) - \esnf_{z \in Q' }\mathcal{M}^{\beta_2} f (z) \right)\dcb{\beta_1}(y) \\
   \leq & \frac{1}{\contb{\beta_1}(Q')} \int_{Q'} \left(\mathcal{M}^{\beta_2} g(y) 
    +  \left( C \|f\|_{\bmob{\beta_2}(\Rn)} \right) \right)\dcb{\beta_1}(y) \\
    \leq &\frac{1}{\contb{\beta_1}(Q')} \int_{Q'} \mathcal{M}^{\beta_2} g(y) \dcb{\beta_1}(y) 
    +  C \|f\|_{\bmob{\beta_2}(\Rn)} \\
    \leq & C \left(\|f\|_{\bmob{\beta_1}(\Rn)}+\|f\|_{\bmob{\beta_2}(\Rn)}\right) \leq C \|f\|_{\bmob{\beta_1}(\Rn)}.
\end{align*}
In the end, we used Lemma~\ref{nesting}, and this completes the proof.
\end{proof}

\begin{proof}[Proof of Theorem~\ref{asthe}]
  Theorem~\ref{asthe} is now a simple consequence of Theorem~\ref{MainThm} by choosing $\beta_1 = \beta$, $\beta_2=n$, and assuming that the support of $f$ is contained in the cube $Q$.
\end{proof}

\section{Vanishing mean oscillation spaces and Fractional maximal function \texorpdfstring{$\fma$}{fma}}\label{Sec: VMO-boundedness}

This section is dedicated to the proofs of Theorem~\ref{VMO-improvement}. In the proof of Theorem~\ref{VMO-improvement}, we require certain key mapping properties of the Riesz potential $I_{\alpha}$. For $\alpha \in (0, n)$ and $f \in L^1_{\mathrm{loc}}(\Rn)$, we recall that the Riesz potential $I_\alpha$ of order $\alpha$ is defined by \begin{align*} 
I_\alpha f(x) = \int_{\Rn} \frac{f(y)}{|x - y|^{n - \alpha}} \dy. 
\end{align*}

The critical Sobolev embedding theorem, formulated in terms of the Riesz potential, states that if $f$ is a function in the Lebesgue space $L^{\frac{n}{\alpha}}(\Rn)$ with compact support, then its Riesz potential $I_\alpha f$ belongs to  $\cbmo(\Rn)$. More precisely, let $0 < \alpha < n$. If $f \in L^{\frac{n}{\alpha}}(\Rn)$ and has compact support, then there exists a constant $C > 0$, depending only on $\alpha$ and $n$, such that 
\begin{align}\label{criticalSobolev} \Vert I_\alpha f \Vert_{\cbmo(\Rn)} \leq C \Vert f \Vert_{L^{\frac{n}{\alpha}}(\Rn)}. 
\end{align}
 For a detailed proof, we refer the reader to \cite{MR3381284}*{Corollary~14.46} and the discussion on \cite{MR3381284}*{p.~446}. In this context, we also mention the original work of Adams \cite{Adams_1975}, where such embeddings are studied within the broader framework of Morrey spaces. Furthermore, for a general trace setting, the critical Sobolev embedding for the Riesz potential in Lorentz spaces is established in the recent work \cite{PD}*{Theorem~1.2}.
 
We require an additional auxiliary result from \cite{Adams_1975}*{Corollary~(i), p.~772}; see also the extension provided in \cite{YC}*{Theorem~1.6}. Suppose $0 < \alpha < n$, and let $f$ be a nonnegative measurable function such that $I_\alpha f \in \cbmo(\Rn)$. Then the following equivalence holds:
\begin{align}\label{Adamsequivalent}
	\Vert I_\alpha f \Vert_{\cbmo(\Rn)} \cong \Vert \fma f \Vert_{L^\infty(\Rn)}.
\end{align}

Now, we give the proof of Theorem \ref{VMO-improvement}.
  
\begin{proof}[Proof of Theorem \ref{VMO-improvement}]
By Lemma~\ref{abs-bmo} and the definition of the fractional maximal function $\fma$, we may assume without loss of generality that $f$ is nonnegative. Let $h \in \Rn$, and consider the translation of $f$, denoted by $\tau_h f$, defined for every $x \in \Rn$ as
\begin{align*}
	(\tau_h f)(x) = f(x + h).
\end{align*}
In particular, by the definition of the fractional maximal function, we observe that it commutes with translations; that is, for every $x \in \Rn$, we have
\begin{align*}
	\fma(\tau_h f)(x) = \tau_h(\fma f)(x).
\end{align*}

Using the above observation and applying the triangle inequality, we obtain
\begin{align}\label{tran}
	\left\Vert \tau_h(\mathcal{M}_\alpha f) - \mathcal{M}_\alpha f \right\Vert_{L^\infty(\Rn)} = \left\Vert \mathcal{M}_\alpha(\tau_h f) - \mathcal{M}_\alpha f \right\Vert_{L^\infty(\Rn)} 
	\leq \left\Vert \mathcal{M}_\alpha(\tau_h f - f) \right\Vert_{L^\infty(\Rn)}.
\end{align}
Since $f$ is supported in $Q_0$ and belongs to $\cbmo(Q_0)$, the John--Nirenberg inequality for $\cbmo(Q_0)$ implies that $f \in L^{\frac{n}{\alpha}}(\Rn)$. Furthermore, applying the John--Nirenberg inequality again, we have $f \in L^p(\Rn)$ for some $1 < p < \frac{n}{\alpha}$. Then, by \cite{MR3381284}*{Theorem~14.37}, it follows that $I_\alpha f \in L^q(\Rn)$ for some $q > 1$, where the exponents satisfy $\frac{1}{q} = \frac{n - \alpha}{np}$. In particular, this implies that $I_\alpha(\tau_h f - f) \in L^1_{\mathrm{loc}}(\Rn)$. 

Then, applying \eqref{Adamsequivalent} and \eqref{criticalSobolev}, we have
\begin{align*}
	\Vert \fma(\tau_h f - f) \Vert_{L^\infty(\Rn)} \leq C \Vert I_\alpha(\tau_h f - f) \Vert_{\cbmo(\Rn)}
	\leq C \Vert \tau_h f - f \Vert_{L^{\frac{n}{\alpha}}(\Rn)},
\end{align*}
where $C$ is a positive constant depending only on $\alpha$ and $n$. Since $f \in L^{\frac{n}{\alpha}}(\Rn)$, it follows that $\tau_h f \to f$ in $L^{\frac{n}{\alpha}}(\Rn)$ as $|h| \to 0$. Thus, we deduce
\begin{align}\label{uni3}
	\lim_{|h| \to 0} \Vert \mathcal{M}_\alpha(\tau_h f - f) \Vert_{L^\infty(\Rn)} \leq C \lim_{|h| \to 0} \Vert \tau_h f - f \Vert_{L^{\frac{n}{\alpha}}(\Rn)} = 0.
\end{align}
Combining \eqref{tran} and \eqref{uni3}, we conclude that
\begin{align*}
	\lim_{|h| \to 0} \left\Vert \tau_h(\fma f) - \fma f \right\Vert_{L^\infty(\Rn)} = 0,
\end{align*}
which implies that $\fma f$ is uniformly continuous on $\Rn$.

Now, by Theorem~\ref{Kline-Gibara-BMO-result}, we have $\fma f \in \cbmo(Q_0)$. Finally, applying \cite{fractionalmaximalpaper}*{Lemma~3.6}, it follows that $\fma f \in \cvmo(Q_0)$, which completes the proof.
\end{proof}

\section{ Vanishing \texorpdfstring{$\beta$}{Beta}-dimensional mean oscillation spaces and fractional maximal function \texorpdfstring{$\fma$}{}}\label{Sec:Beta-VMO}
The primary objective of this section is to analyze the behavior of the fractional maximal function within the framework of the vanishing $\beta$-dimensional mean oscillation space. In particular, we focus on the proof of Theorem~\ref{thm-VMO-main}. The proof relies on decomposing the fractional maximal function $\fma$ into its local and global components. We then provide a detailed estimation of the mean oscillation corresponding to each component to establish the desired result.

We begin by presenting several key lemmas that will be instrumental in the proof of Theorem~\ref{thm-VMO-main}, the main result of this section.
\begin{lemma}\label{tri_const_cont}
	Suppose $0<\beta\leq n$. Let $Q$ be any finite cube in $\Rn$, and let $f \in L^1(\cont, Q)$ be a nonnegative function. Then, for any constant $k \in \mathbb{R}$, the following inequality holds
	\begin{align*}
		\frac{1}{\cont(Q)}\int_{Q} |f - f_{Q,\beta}| \, \dc \leq \frac{2}{\cont(Q)}\int_{Q} |f - k| \, \dc. 
	\end{align*}
\end{lemma}
\begin{proof}
		Let us assume that $k \geq 0$. Now, using the sublinearity of the integral and Lemma~\ref{tri_cont}, we obtain
		\begin{align*}
			&	\frac{1}{\cont(Q)}\int_{Q} |f-\avgd{Q}|\dc \\&\leq \frac{1}{\cont(Q)}\int_{Q^\prime} |f-k|\dc + \frac{1}{\cont(Q)}\int_{Q} |k-\avgd{Q}|\dc	\\& = \frac{1}{\cont(Q)}\int_{Q} |f-k|\dc + |k-\avgd{Q}|\\& =  \frac{1}{\cont(Q)}\int_{Q} |f-k|\dc + \left|\frac{1}{\cont(Q)}\int_{Q} |k|\dc-\frac{1}{\cont(Q)}\int_{Q}|f|\dc \right|\\& \leq \frac{2}{\cont(Q)}\int_{Q} |f-k|\dc.
		\end{align*}
		In particular, for $k = 0$, we have
		\begin{align*}
			\frac{1}{\cont(Q)}\int_{Q} |f-\avgd{Q}|\dc \leq  \frac{2}{\cont(Q)}\int_{Q} |f|\dc.
		\end{align*}
	Now, for other values of $k < 0$, we have $|f| \leq |f - k|$ pointwise. Applying this inequality together with the result above and the monotonicity of the integral, we conclude the proof.
	\end{proof}

\begin{remark}\label{c-avg}
	Suppose $0 < \beta \leq n$, and let $Q$ be any finite cube in $\Rn$. If $f \in L^1(\cont, Q)$ is a nonnegative function, then, by Lemma~\ref{tri_const_cont}, we have
	\begin{align*}
		\O_{\beta}(f, Q) \cong \frac{1}{\cont(Q)}\int_{Q} |f - \avgd{Q}| \, \dc.
	\end{align*}
\end{remark}

\begin{lemma}\label{vmo-comp}
	Let $0 < \beta_1 \leq \beta_2 \leq n$. Then, for any cube $Q$, the following inequality holds:
	\begin{align*}
		\O_{\beta_2}(f, Q) \leq C_0 \, \O_{\beta_1}(f, Q),
	\end{align*}
	where $C_0$ is a positive universal constant independent of $Q$. Consequently, $\vmob(Q_0) \subset \text{VMO}^{\beta_2}(Q_0)$, where $Q_0$ is any finite open cube in $\Rn$.
\end{lemma}

\begin{proof}
	The argument follows by mimicking the proof of Lemma~\ref{compare-beta} and applying the sublinearity of the Choquet integral.   
\end{proof}
\begin{lemma}\label{osc-add}
	Suppose $0 < \beta \leq n$, and let $Q_0$ be a finite cube in $\Rn$. Assume $f, g \in L^1_{\mathrm{loc}}(\cont, Q_0)$, and define $h(x) = \max\{f(x), g(x)\}$ for $x \in Q_0$. Then the following hold
	\begin{enumerate}
		\item $h \in L^1_{\mathrm{loc}}(\cont, Q_0)$;
		\item for every finite cube $Q \subset Q_0$, we have
		\begin{align*}
			\O_{\beta}(h, Q) \leq \O_{\beta}(f, Q) + \O_{\beta}(g, Q).
		\end{align*}
	\end{enumerate}
\end{lemma}
\begin{proof}
		To prove part (1), observe that since $|h| \leq |f| + |g|$ and the Choquet integral with respect to $\cont$ is sublinear, it follows $h \in L^1_{\mathrm{loc}}(\cont, Q_0)$.
		
		For part (2), fix a constant $c \in \mathbb{R}$. For any cube $Q \subset Q_0$, using the subadditivity of the absolute value and the sublinearity of the Choquet integral, we have  
		\begin{align*}
			&	\frac{1}{\cont(Q)} \int_{Q}|h-c|\dc \\& \leq \frac{1}{\cont(Q)} \int_{0}^\infty \cont \left(\{x\in Q\,:\, |f(x)-c|>t , \, f(x)\geq g(x)\}\right)\dt \\& + \frac{1}{\cont(Q)} \int_{0}^\infty \cont \left(\{x\in Q\,:\, |g(x)-c|>t , \, f(x)< g(x)\}\right)\dt \\& \leq \frac{1}{\cont(Q)} \int_{0}^\infty \cont \left(\{x\in Q\,:\, |f(x)-c|>t \}\right)\dt \\& + \frac{1}{\cont(Q)} \int_{0}^\infty \cont \left(\{x\in Q\,:\, |g(x)-c|>t\}\right)\dt \\& = \frac{1}{\cont(Q)} \int_{Q}|f-c|\dc+ \frac{1}{\cont(Q)} \int_{Q}|g-c|\dc.
		\end{align*}
	In the end, taking the infimum over all constants $c \in \mathbb{R}$, we conclude 
		\begin{align*}
			\O_{\beta}(h,Q)\leq \O_{\beta}(f,Q)+\O_{\beta}(g,Q).
		\end{align*}
	\end{proof}

   The next lemma ensures that it is sufficient to consider nonnegative functions in our analysis.
\begin{lemma}\label{lem-abs-VMO}
	Suppose $0 < \beta \leq n$, and let $Q_0 \subset \Rn$ be a finite cube. If $f \in \vmo(Q_0)$, then $|f| \in \vmo(Q_0)$.
\end{lemma}
 	 
	\begin{proof}
		Since $f \in \vmo(Q_0)$, it follows by definition that $f \in \bmo(Q_0)$. Hence, by Lemma~\ref{abs-bmo}, we have $|f| \in \bmo(Q_0)$. Moreover, the vanishing $\beta$-dimensional mean oscillation condition implies that
		\begin{align}\label{eq-VMO-lim}
			\lim_{r \to 0} \omega_{\beta}(f, r) = 0.
		\end{align}
		Let us now fix a cube $Q \subset Q_0$. For $x \in Q$, using the same calculation as in Lemma~\ref{abs-bmo}, we obtain
		\begin{align*}
			\left| |f(x)| - |f|_{Q, \beta} \right| 
			\leq |f(x) - c_Q| + \frac{1}{\cont(Q)} \int_Q |f(y) - c_Q| \, \dc 
			= |f(x) - c_Q| + \O_{\beta}(f, Q).
		\end{align*}
		
		Now, taking the average over $Q$, we have
		\begin{align*}
			\frac{1}{\cont(Q)} \int_Q \left| |f(x)| - |f|_{Q, \beta} \right| \, \dc 
			\leq \frac{1}{\cont(Q)} \int_Q |f(x) - c_Q| \, \dc + \O_{\beta}(f, Q) 
			= 2 \, \O_{\beta}(f, Q).
		\end{align*}
		
		Moreover, we observe that
		\begin{align}\label{eq-infy-Vmo}
			\inf_{c \in \mathbb{R}} \frac{1}{\cont(Q)} \int_Q \left| |f(x)| - c \right| \, \dc 
			\leq \frac{1}{\cont(Q)} \int_Q \left| |f(x)| - |f|_{Q, \beta} \right| \, \dc 
			\leq 2 \, \O_{\beta}(f, Q).
		\end{align}
	
	Taking the supremum over all cubes $Q \subset Q_0$ with $\ell(Q) < r$ on both sides of \eqref{eq-infy-Vmo}, we obtain
	\begin{align*}
		\omega_{\beta}(|f|, r) \leq 2 \, \omega_{\beta}(f, r).
	\end{align*}
	Now, taking the limit as $r \to 0$ and using \eqref{eq-VMO-lim}, we conclude that
	\[
	\lim_{r \to 0} \omega_{\beta}(|f|, r) = 0.
	\]
	Hence, $|f| \in \vmo(Q_0)$. This completes the proof of the lemma. 
	\end{proof} 

We continue by presenting a series of lemmas that highlight some useful observations concerning cubes.
\begin{lemma}\label{interval}
	 Let $I_1$ and $I_2$ be intervals in $\mathbb{R}$ such that $I_1 \subseteq I_2$. Then, for every $\delta > 0$ satisfying $ \ell(I_1) \leq \delta \leq \ell(I_2)$, there exists an interval $I'$ such that: \begin{enumerate} \item $I_1 \subseteq I' \subseteq I_2$; \item $\ell(I') = \delta$. \end{enumerate}
 \end{lemma}
\begin{proof} The result follows directly from the properties of intervals in the real line.
 \end{proof}

\begin{lemma}\label{trivialfact1} Let $Q$ be a finite cube in $\mathbb{R}^n$. If $Q_1$ and $Q_2$ are cubes contained in $Q$ such that $Q_1 \cap Q_2 \neq \emptyset$, then there exists a cube $Q'$ satisfying: \begin{enumerate} \item $Q' \supseteq Q_1 \cup Q_2$; \item $Q' \subseteq Q$; \item $\ell(Q') \leq \ell(Q_1) + \ell(Q_2)$. \end{enumerate} \end{lemma}

\begin{proof} We represent the cubes $Q$, $Q_1$, and $Q_2$ as Cartesian products of intervals in $\mathbb{R}$: \begin{align*} Q = \prod_{i=1}^n I_i, \quad Q_1 = \prod_{i=1}^n I_i^1, \quad Q_2 = \prod_{i=1}^n I_i^2, \end{align*} where $I_i$, $I_i^1$, and $I_i^2$ are intervals in $\mathbb{R}$. Since $Q_1 \cap Q_2 \neq \emptyset$, the intervals $I_i^1$ and $I_i^2$ intersect non-trivially in each coordinate. Thus, the union $I_i^1 \cup I_i^2$ is an interval contained in $I_i$. Define \begin{align}\label{deltalength} \delta = \max_{1 \leq i \leq n} \ell(I_i^1 \cup I_i^2). \end{align} Observe that for all $1 \leq i \leq n$, \begin{align*} \ell(I_i^1 \cup I_i^2) \leq \ell(I_i^1) + \ell(I_i^2) = \ell(Q_1) + \ell(Q_2), \end{align*} and therefore, \begin{align*} \delta \leq \ell(Q_1) + \ell(Q_2). \end{align*}
	
	Moreover, since $I_i^1 \cup I_i^2 \subseteq I_i$, it follows that $\ell(I_i^1 \cup I_i^2) \leq \delta \leq \ell(I_i)$ for each $i$. By applying Lemma~\ref{interval} to each coordinate $i$, with intervals $I_i^1 \cup I_i^2 \subseteq I_i$ and the length $\delta$ as defined in (\ref{deltalength}), we obtain intervals $I_i'$ satisfying: \begin{enumerate} \item $I_i^1 \cup I_i^2 \subseteq I_i' \subseteq I_i$; \item $\ell(I_i') = \delta$. \end{enumerate} Define $Q' = \prod_{i=1}^n I_i'$. Then $Q'$ is a cube of side length $\delta$, satisfying: \begin{enumerate} \item $Q' \supseteq Q_1 \cup Q_2$, since $I_i' \supseteq I_i^1 \cup I_i^2$ for all $i$; \item $Q' \subseteq Q$, because $I_i' \subseteq I_i$ for all $i$; \item $\ell(Q') = \delta \leq \ell(Q_1) + \ell(Q_2)$. \end{enumerate} 
    This completes the proof of the lemma. 
    \end{proof}

We have established all the necessary preliminary results and now proceed to prove Theorem~\ref{thm-VMO-main}. For this purpose, we decompose the fractional maximal function $\fma$ into two components: a local part and a global part, under the assumptions $0 \leq \alpha < n$. We fix a scaling parameter $\kappa > 0$. Then, we have
\begin{align*}
		\fma f = \max\{\fmal  f \, , \,  \fmag f\},
\end{align*}
where we define the local part of $\fma f$ as 
	\begin{align*}
		\fmal  f(x) :=\sup_{\{  x\in Q^\prime, \, \ell(Q^\prime)<\kappa\}} \frac{1}{\ell(Q^\prime)^{n-\alpha}}\int_{Q^\prime}|f(y)|\dy,
	\end{align*}
	and the global part  of $\fma f$ as 
	\begin{align*}
		\fmag  f(x) :=\sup_{\{ x\in Q^\prime, \, \ell(Q^\prime)\geq \kappa\}} \frac{1}{\ell(Q^\prime)^{n-\alpha}}\int_{Q^\prime}|f(y)|\dy.
	\end{align*}
In the above definitions, the supremum is taken over all cubes $Q \subseteq \mathbb{R}^n$ containing the point $x$, with side length constrained appropriately in each case.
\begin{remark}
Suppose $0 \leq \alpha < n$, and let $Q_0$ be a finite cube in $\mathbb{R}^n$. Assume that $\supp f \subset Q_0$. Then, by applying Lemma~\ref{diffBMO} and using arguments similar to those in \eqref{ywl-mf}, we observe that for any $\kappa > 0$ and for all $x \in Q_0$, one has
\begin{align}\label{equialocmaxima}
    \fmal  f (x) = \sup_{\{ x\in Q, \ell(Q)\leq \kappa, Q\subseteq Q_0\}} \frac{1}{\ell(Q)^{n-\alpha} } \int_Q |f(y)| \dy,
\end{align} 
and 
\begin{align}\label{equiagolbmaxima}
     \nonumber\fmag  f(x)=&\sup_{\{ x\in Q, \kappa<\ell(Q)\leq \ell(Q_0)\}} \frac{1}{\ell(Q)^{n-\alpha}} \int_Q |f(y)| \dy \\
     =&\sup_{\{ x\in Q, \kappa<\ell(Q)\leq \ell(Q_0), Q\subseteq Q_0\}} \frac{1}{\ell(Q)^{n-\alpha}} \int_Q |f(y)| \dy.
\end{align}
\end{remark}

Next, we will establish two auxiliary lemmas essential in proving our main result. 
\begin{lemma}\label{loc-osc}
	Let $0 < \beta \leq n$, $\alpha \in [0, \beta)$, and let $Q_0 \subset \mathbb{R}^n$ be a finite open cube. Suppose $f$ is a nonnegative function supported in $Q_0$, and assume that $f \in L^1_{{\rm loc}}(\cont,Q_0)$. Then for any cube $\tilde{Q} \subset Q_0$ with side length $\ell(\tilde{Q}) = r$ and $\lambda\geq 1$,  there exists a universal constant $C$ such that the following holds
	\begin{align*}
		\O_{\beta}(\maxl f ,\tilde{Q}) \leq C \lambda^{\frac{\beta}{p}} \ell(Q_0)^\alpha \omega_{\beta}(f, 3\lambda r),
	\end{align*}
    where $1<p<\frac{\beta}{\alpha}$. 
\end{lemma}
\begin{proof}
	Let $\tilde{Q} \subset Q_0$ be a cube with $\ell(\tilde{Q}) = r$. Then, by Lemma~\ref{diffBMO} with $\gamma=0$, there exists a cube $P_{3\lambda\tilde{Q}}$ associated with the dilated cube $3\lambda\tilde{Q}$ such that the following properties hold:
\begin{enumerate}[label=(\roman*)]
\item $P_{3\lambda\tilde{Q}} \subseteq Q_0 \cap  6\lambda\tilde{Q}$.  
\item $P_{3\lambda\tilde{Q}} \supseteq Q_0 \cap  3\lambda \tilde{Q}$.  
\item The side length of $P_Q$ satisfies
\begin{align*}
    \ell(P_{3\lambda\tilde{Q}}) = \begin{cases}
   \ell(3\lambda\tilde{Q}), & \text{if } \ell(3\lambda\tilde{Q}) < \ell(Q_0),\\
   \ell(Q_0), & \text{if } \ell(3\lambda\tilde{Q}) \geq \ell(Q_0).
   \end{cases}
\end{align*}

\item Fix any $p \in \big( \frac{\beta}{n}, \frac{\beta}{\alpha} \big)$. Then the following inequality holds:
\begin{align*}
    \left( \frac{1}{\contb{\beta}(3\lambda\tilde{Q})} \int_{3\lambda\tilde{Q} \cap  Q_0} |f - \avgdb{P_{3\lambda\tilde{Q}}}{\beta}|^p \dcb{\beta} \right)^{\frac{1}{p}}
   \leq D \left( \frac{1}{\contb{\beta}(P_{3\lambda\tilde{Q}})} \int_{P_{3\lambda\tilde{Q}}} |f - \avgdb{P_{3\lambda\tilde{Q}}}{\beta}|^p \dcb{\beta}  \right)^{\frac{1}{p}}.
\end{align*}
\end{enumerate}

Next, by the definition of oscillation and the definition of $\maxl f$, we obtain the following
\begin{align}\label{initial-cal}
			\O_{\beta}(\maxl f ,\tilde{Q})&=\inf_{c \in \mathbb{R}}\frac{1}{\contb{\beta}(\tilde{Q})}\int_{\tilde{Q}}|\maxl f-c|\dcb{\beta} \\& \leq \frac{1}{\contb{\beta}(\tilde{Q})}\int_{\tilde{Q}}\left|\maxl f-(\lambda r)^\alpha \avgd{P_{3\lambda\tilde{Q}}}\right|\dcb{\beta}\nonumber \\&  =\frac{1}{\contb{\beta}(\tilde{Q})}\int_{\tilde{Q}}\left|\maxl f-\maxl (\avgd{P_{3\lambda\tilde{Q}}})\right|\dcb{\beta}. \nonumber
\end{align}
Now, take any $x \in \tilde{Q}$. We observe that for any cube $Q'$ with $x \in Q'$ and $\ell(Q') < \lambda r$, the condition $x \in \tilde{Q}$ and $\ell(\tilde{Q}) = r$ implies that $Q' \subset 3\lambda\tilde{Q}$. Therefore, any such cube $Q'$ containing $x$ lies within $3\lambda \tilde{Q}$. Furthermore, since $x \in \tilde{Q} \subset Q_0$ and $\supp f \subset Q_0$, we can apply equation~\eqref{equialocmaxima}, and hence we can write
\begin{align}\label{obs-1}
    \nonumber\maxl (f\chi_{3\lambda\tilde{Q}\cap Q_0})(x) &=\sup_{x\in Q^\prime, \, \ell(Q^\prime)<\lambda r} \frac{1}{\ell(Q')^{n-\alpha}}\int_{Q^\prime}|f\chi_{3\lambda \tilde{Q}\cap Q_0}(y)|\dy\nonumber\\
    &=\sup_{x\in Q^\prime, \, \ell(Q^\prime)<\lambda r}\frac{1}{\ell(Q')^{n-\alpha}}\int_{Q^\prime\cap 3\lambda\tilde{Q}\cap  Q_0}|f(y)|\dy \nonumber\\
    &=\sup_{x\in Q^\prime, \, \ell(Q^\prime)<\lambda r} \frac{1}{\ell(Q')^{n-\alpha}}\int_{Q^\prime\cap  Q_0}|f(y)|\dy \nonumber\\&=\sup_{x\in Q^\prime, \, \ell(Q^\prime)<\lambda r,\, Q^\prime\subset Q_0 }\frac{1}{\ell(Q')^{n-\alpha}}\int_{Q^\prime}|f(y)|\dy\nonumber\\&=\maxl f(x).
  		\end{align}
  	
Since $3\lambda\tilde{Q} \cap Q_0 \subset Q_0$, we can apply the equation~\eqref{equialocmaxima} and deduce that
\begin{align}\label{obs-2}
	\nonumber&\maxl \left(\avgdb{P_{3\lambda\tilde{Q}}}{\beta}\chi_{3\lambda\tilde{Q}\cap Q_0}\right)(x)\\ &=\avgdb{P_{3\lambda\tilde{Q}}}{\beta} \sup_{x\in Q^\prime, \, \ell(Q^\prime)<\lambda r, Q'\subset Q_0} \frac{1}{\ell(Q')^{n-\alpha}}\int_{Q^\prime}|\chi_{3\lambda \tilde{Q}\cap Q_0}(y)|\dy \nonumber\\&=\avgdb{P_{3\lambda\tilde{Q}}}{\beta} \sup_{x\in Q^\prime, \, \ell(Q^\prime)<\lambda r, Q'\subset Q_0} \frac{1}{\ell(Q')^{n-\alpha}}\int_{Q^\prime}\dy \nonumber\\& =\maxl (\avgdb{P_{3\lambda\tilde{Q}}}{\beta})(x),
\end{align}
holds for any $x \in \tilde{Q}.$

In the second equality, we used the fact that the points outside the cube $3\lambda\tilde{Q}$ cannot be covered by any cube of side length at most $\lambda r$ that contains $x$.
        
        Hence, using the observations from \eqref{obs-1}, \eqref{obs-2}, and \eqref{initial-cal}, we get
		\begin{align*}
			\O_{\beta}(\maxl f ,\tilde{Q}) & \leq \frac{1}{\contb{\beta}(\tilde{Q})}\int_{\tilde{Q}}\left|\maxl f-\maxl (\avgdb{P_{3\lambda\tilde{Q}}}{\beta})\right|\dcb{\beta} \\& = \frac{1}{\contb{\beta}(\tilde{Q})}\int_{\tilde{Q}}\left|\maxl (f\chi_{3\lambda\tilde{Q}\cap Q_0})-\maxl (\avgdb{P_{3\lambda\tilde{Q}}}{\beta}\chi_{3\lambda \tilde{Q}\cap Q_0})\right|\dcb{\beta}\\&\leq \frac{1}{\contb{\beta}(\tilde{Q})}\int_{\tilde{Q}} \maxl [(f-\avgdb{P_{3\lambda\tilde{Q}}}{\beta})\chi_{3\lambda \tilde{Q}\cap Q_0}]\dcb{\beta}\\
            & \leq  \frac{1}{\contb{\beta}(\tilde{Q})} \int_{\tilde{Q}} \fma g \dcb{\beta}\\ & \lesssim \contb{\beta}(\tilde{Q})^{\frac{p\alpha-\beta}{\beta p}} ||g||_{L^p(\contb{\beta}, \Rn)}, 
		\end{align*}
        where $g=(f-\avgdb{P_{3\lambda\tilde{Q}}}{\beta})\chi_{3\lambda\tilde{Q}\cap Q_0}$. In the last inequality, we used H\"older's inequality and  \cite{Herjulehto-Petteri_2024}*{Theorem~4.15} for $\delta=\beta$, $k=\alpha$ with $1<p<\frac{\beta}{\alpha}$. Notice that the centers of the cubes $3\lambda\tilde{Q}$ and $\lambda\tilde {Q}$ coincide with the center of $\tilde{Q}$, and since $\tilde{Q} \subset Q_0$, this center is contained within $Q_0$. Now, rewriting the expression for $g$, and using (iv), along with the doubling condition and equation~\eqref{equiv-cont-cube}, we obtain 
		\begin{align*}
			\O_{\beta}(\maxl f ,\tilde{Q})  & \lesssim \lambda^{\frac{\beta}{p}}\contb{\beta}(\tilde{Q})^{\frac{\alpha}{\beta }}\left( \frac{1}{\cont(3\lambda\tilde{Q})}\int_{3\lambda \tilde{Q}\cap Q_0}|f-\avgd{P_{3\lambda\tilde{Q}}}|^p\dcb{\beta} \right)^{1/p} \\ & \lesssim  \lambda^{\frac{\beta}{p}} \ell(Q_0)^\alpha \left( \frac{1}{\contb{\beta}(3\lambda\tilde{Q})}\int_{3\lambda \tilde{Q}\cap Q_0}|f-\avgdb{P_{3\lambda\tilde{Q}}}{\beta}|^p\dcb{\beta} \right)^{1/p}\\& \lesssim \lambda^{\frac{\beta}{p}}\ell(Q_0)^\alpha  \left( \frac{1}{\contb{\beta}(P_{3\lambda\tilde{Q}})} \int_{P_{3\lambda\tilde{Q}}} |f - \avgdb{P_{3\lambda\tilde{Q}}}{\beta}|^p \dcb{\beta}  \right)^{1/p}.
		\end{align*}

Therefore, finally, using Lemma~\ref{tri_const_cont}, the information from (i), namely $P_{3\lambda\tilde{Q}} \subset Q_0$, along with the application of the John-Nirenberg inequality \cite{Chen-Spector}*{Corollary~1.5}, and utilizing the previous results, we deduce   
		\begin{align*}
			\O_{\beta}(\maxl f ,\tilde{Q})&\lesssim \lambda^{\frac{\beta}{p}} \ell(Q_0)^\alpha \|f\|_{\bmobp{\beta}{p}(P_{3\lambda \tilde{Q}}) }\leq C \lambda^{\frac{\beta}{p}} \ell(Q_0)^\alpha \|f\|_{\bmob{\beta}( P_{3\lambda\tilde{Q}})} \leq C \lambda^{\frac{\beta}{p}} \ell(Q_0)^\alpha \omega_{\beta}(f, 3\lambda r),
		\end{align*}
	where $C$ depends on universal parameters. In the final step, we use the information from (iii), namely that $\ell(P_{3\lambda\tilde{Q}}) \leq \ell(3\lambda\tilde{Q}) = 3\lambda r$. This completes the proof of the lemma.
    \end{proof}

\begin{lemma}\label{glob-osc}
   Let $0 < \beta \leq n$, $\alpha \in [0, \beta)$, and let $Q_0 \subset \mathbb{R}^n$ be a finite open cube. 
Assume \(f\in\bmo(Q_0)\) is non-negative and supported in \(Q_0\).
Then there exists a dimensional constant \(C_n>1\) such that for every cube
\(\tilde Q\subset Q_0\) with side length \(\ell(\tilde Q)=r\) and for every
\(\lambda>C_n e\), one can find a constant \(C>0\) (independent of $\lambda$ and \(r\)) for which
\begin{align*}
		\O_{\beta}(\maxgg f, \tilde{Q}) \leq C \ell(Q_0)^\alpha \lambda^{-1} \left(1+\log \lambda\right)\|f\|_{\bmo(Q_0)}.
	\end{align*}
\end{lemma}
\begin{proof}
Let us take $x,y \in \tilde{Q}$. Without loss of generality, we may assume that
\begin{align}\label{maxi-comp}
		 \maxgg f(y)<\maxgg f(x).
\end{align}
Therefore, by the inequality in \eqref{maxi-comp}, to estimate the mean oscillation of the global part, we choose a cube $Q_x$ with $\ell(Q_x) \geq \lambda r$ such that $x \in Q_x$ and  
		\begin{align}\label{gml}
			\maxgg f(y) <  \ell(Q_x)^{\alpha}\avgl{Q_x}.
		\end{align}
	Additionally, by (\ref{equiagolbmaxima}), we may assume $Q_x \subset Q_0$. Now, applying Lemma~\ref{trivialfact1}, we can choose a cube $Q'$ such that $Q' \subset Q_0$, $\tilde{Q} \cup Q_x \subset Q'$, and
		\begin{align}\label{vmog-1}
			 \ell(Q') \leq \ell(Q_x) + \ell(\tilde{Q}).
		  \end{align}

 Since $Q_x \subset Q'$, it follows that $\lambda r \leq \ell(Q_x) \leq \ell(Q')$. Moreover, as $y \in \tilde{Q} \subset Q'$, we obtain that $\maxgg f(y) \geq \ell(Q')^{\alpha} \avgl{Q'}$. Now, consider $\Gamma:=Q'\setminus Q_x$  and set $$f_\Gamma : = \frac{1}{\Gamma} \int_\Gamma f \dx .  $$ A straightforward computation then gives
		\begin{align}\label{vmog-2}
			\nonumber \ell(Q_x)^{\alpha}\avgl{Q_x}-\maxgg f(y)&\leq \ell(Q_x)^{\alpha}\avgl{Q_x}-\ell(Q')^{\alpha}\avgl{Q'}\\&\nonumber \leq \ell(Q_x)^{\alpha}\left(\avgl{Q_x}-\avgl{Q'}\right) \\& \nonumber \leq \ell(Q_0)^{\alpha}\frac{|\Gamma|}{\ell(Q')^n}\left(\avgl{Q_x}-\avgl{\Gamma}\right)\\& \nonumber \leq \ell(Q_0)^{\alpha}\frac{|\Gamma|}{\ell(Q')^n}\left(|\avgl{Q_x}-\avgl{Q'}|+|\avgl{Q'}-\avgl{\Gamma}|\right)\\
             &\leq \ell(Q_0)^{\alpha}\left(J_1+J_2\right),
		\end{align}	
        where
        \begin{align*}
            J_1:=\frac{|\Gamma|}{\ell(Q_x)^n} \frac{1}{\ell(Q')^n}\int_{Q'}|f-\avgl{Q'}| \dy \quad \text{ and } \quad J_2:=\frac{1}{\ell(Q')^n} \int_{\Gamma}|f-\avgl{Q'}| \dy.
        \end{align*}
        Let us begin with the estimate of $J_1$. We observe
        \begin{align}\label{obs}
            \frac{|\Gamma|}{\ell(Q_x)^n}=\frac{\ell(Q')^n-\ell(Q_x)^n}{\ell(Q_x)^n}\leq \frac{\left(\ell(Q_x)+\ell(\tilde{Q})\right)^n-\ell(Q_x)^n}{\ell(Q_x)^n}\leq C_n \frac{\ell(\tilde{Q})}{\ell(Q_x)} \leq C_n\lambda^{-1},
        \end{align}
        where $C_n$ is a positive constant that depends on the dimension $n$. Now using this, Remark~\ref{c-avg}, and Lemma~\ref{nesting}, we have 
        \begin{align}\label{eye1}
          J_1\lesssim \lambda^{-1}\O_n(f,Q')\lesssim \lambda^{-1}\|f\|_{\bmo(Q_0)}.  
        \end{align}

        Next, applying \cite{xianbao}*{Lemma~3.2} with $A=\Gamma$ and $Q=Q'$, we deduce
        \begin{align*}
            \frac{1}{|\Gamma|} \int_{\Gamma}|f-\avgl{Q'}| \dy \lesssim \|f\|_{\bmo(Q_0)} \left(1+\log\frac{\ell(Q')^n}{|\Gamma|} \right).
        \end{align*}
        Hence, using the estimate $\ell (Q_x) \leq \ell(Q')\leq 2 \ell(Q_x)$ and \eqref{obs} in $J_2$, we deduce
\begin{align}\label{eye2}
    \nonumber J_2 & \lesssim \|f\|_{\bmo(Q_0)} \frac{|\Gamma|}{\ell(Q')^n}\left(1+\log\frac{\ell(Q')^n}{|\Gamma|} \right) \\
 \nonumber   & =  \|f\|_{\bmo(Q_0)} \frac{|\Gamma|}{\ell(Q')^n}\left(1+ \log \frac{\ell(Q')^n}{\ell(Q_x)^n} 
 +\log\frac{\ell(Q_x)^n}{|\Gamma|} \right)
    \\ &\nonumber \leq \|f\|_{\bmo(Q_0)} \frac{|\Gamma|}{\ell(Q_x)^n}\left( 1+ n \log 2+\log\frac{\ell(Q_x)^n}{|\Gamma|} \right) \\&\nonumber=  \|f\|_{\bmo(Q_0)} \left((1+ n \log 2)\frac{|\Gamma|}{\ell(Q_x)^n}-\frac{|\Gamma|}{\ell(Q_x)^n}\log\frac{|\Gamma|}{\ell(Q_x)^n} \right) \\& \nonumber \leq \|f\|_{\bmo(Q_0)} \left((1+ n \log 2)C_n\lambda^{-1}-C_n\lambda^{-1}\log(C_n\lambda^{-1})\right)\\& \lesssim \|f\|_{\bmo(Q_0)} \left(\lambda^{-1}+\lambda^{-1}\log\lambda \right),
\end{align}
where we have used that the function $-t\log t$ is increasing whenever $t<e^{-1}$  and the assumption $ \lambda>C_ne$.  Finally, combining \eqref{eye1} and \eqref{eye2} in \eqref{vmog-2}, we have
\begin{align*}
	\ell(Q_x)^{\alpha}\avgl{Q_x}-\maxgg f(y) \lesssim \ell(Q_0)^\alpha\|f\|_{\bmo(Q_0)} \lambda^{-1} \left(1+\log \lambda\right).
\end{align*}
Taking the supremum over all such cubes $Q_x$, we conclude
\begin{align*}
	|\maxgg f(x)-\maxgg f(y)|\lesssim \ell(Q_0)^\alpha\|f\|_{\bmo(Q_0)} \lambda^{-1} \left(1+\log \lambda\right) \text{ for all } x, \, y \in \tilde{Q}.
\end{align*}
This yields that on any such cube $\tilde{Q}$, i.e., with $\ell(\tilde{Q}) \leq \tilde{r}$, we have
\begin{align*} 
    \O_{\beta}(\maxgg f, \tilde{Q})&\leq \frac{1}{\contb{\beta}(\tilde{Q})^{2}}\int_{\tilde{Q}}\int_{\tilde{Q}}|\maxgg f(x)-\maxgg f(y)|\dcb{\beta}(x)\dcb{\beta}(y)\\& \lesssim \ell(Q_0)^\alpha \lambda^{-1} \left(1+\log \lambda\right)\|f\|_{\bmo(Q_0)}.
\end{align*} 
This completes the proof of the lemma.
\end{proof}

Now, using the two earlier results, we will show that if $f \in \vmo(Q_0)$, then the oscillations of both the local and global components of $\fma$ can be made arbitrarily small. We will then proceed to demonstrate that the same property holds for the original uncentered fractional maximal function, which we now set out to prove.

\begin{proof}[Proof of the Theorem~\ref{thm-VMO-main}]
	If $f \in \vmo(Q_0)$, then by Lemma~\ref{lem-abs-VMO}, it follows that $|f| \in \vmo(Q_0)$. Moreover, since $\fma f = \fma |f|$, we may assume $f$ is a nonnegative function without loss of generality.
	
	By Theorem~\ref{asthe}, the $\bmob{\beta}(Q_0)$ norm of $\mathcal{M} f$ is finite and by Theorem~\ref{quasi} with $\beta_2=n$ and $\beta_1=\beta$, we have $\mathcal{M} f$ is $\contb{\beta}$-quasicontinuous on $Q_0$. Again by Theorem~\ref{MainThm2} and Lemma~\ref{nesting}, the $\bmob{\beta}(Q_0)$ norm of $\fma f$ is finite whenever $0<\alpha<\beta$. Furthermore, by Theorem~\ref{VMO-improvement}, the function $\fma f$ is uniformly continuous and so $\contb{\beta}$-quasicontinuous on $Q_0$.
	
	Fix a finite cube $\tilde{Q} \subset Q_0$ with $\ell(\tilde{Q}) = r$, and for any $\lambda>C_ne$ (where \(C_n\) is the dimensional constant from
Lemma~\ref{glob-osc}), decompose the maximal function $\fma f$ into a global part $\maxgg f$ and a local part $\maxl f$. Applying Lemma~\ref{osc-add}, we obtain the estimate 
    \begin{align}\label{mf-vmo-eq-0} 
    \O_{\beta}(\fma f, \tilde{Q}) \leq \O_{\beta}(\maxgg f, \tilde{Q}) + \O_{\beta}(\maxl f, \tilde{Q}). 
    \end{align}
	
	Let $\epsilon > 0$ be arbitrary. Using Lemma~\ref{glob-osc}, and using the fact that $\lambda^{-1}\rightarrow 0$ and $\lambda^{-1}\log \lambda\rightarrow 0$ for $\lambda\rightarrow \infty$, we can choose large $\lambda>C_ne$, such that  
    \begin{align}\label{mf-vmo-eq-2} 
        \O_{\beta}(\maxgg f, \tilde{Q}) \leq C \ell(Q_0)^\alpha \lambda^{-1} \left(1+\log \lambda\right)\|f\|_{\bmo(Q_0)}<\frac{\epsilon}{2}.
    \end{align}
    Moreover, we can choose the above large $\lambda$ independent of $\tilde{Q}$. 
    
    Next, we will estimate the local part. Since $f \in \vmo(Q_0)$, this ensures the existence of a constant $r_\epsilon > 0$ such that for any $r \leq r_\epsilon$, there holds
    \begin{align}\label{mf-vmo-eq-1}
    \omega_{\beta}(f, 3\lambda r) < \frac{\epsilon}{2C \lambda^{\frac{\beta}{p}} \ell(Q_0)^\alpha },
    \end{align}
    where $C$ and $p$ are the same as in Lemma~\ref{loc-osc}. Hence, applying Lemma \ref{loc-osc} once more, we conclude that whenever
\(\ell(\tilde Q)=r\le r_\varepsilon\),
\begin{align}\label{mf-vmo-eq-3}
    \O_{\beta}(\maxl f, \tilde{Q}) \leq C \lambda^{\frac{\beta}{p}} \ell(Q_0)^\alpha \omega_{\beta}(f, 3\lambda r) < \frac{\epsilon}{2}.
    \end{align}
	
	Combining \eqref{mf-vmo-eq-3} and \eqref{mf-vmo-eq-2} with \eqref{mf-vmo-eq-0}, we have that for all $\epsilon>0$ there exists $r_\epsilon>0$ such that for any cubes $\tilde{Q}$ whenever $\ell(\tilde{Q}) \leq r_\epsilon$, we get \begin{align*} \O_{\beta}(\fma f, \tilde{Q}) < \epsilon. \end{align*}
	
	Taking the supremum over all such cubes, we obtain \begin{align*} \omega_{\beta}(\fma f, r) < \epsilon \quad \text{for all } r \leq r_\epsilon. \end{align*}
	
	Therefore, by the arbitrariness of $\epsilon > 0$, we conclude that $\fma f \in \vmo(Q_0)$. This completes the proof of Theorem \ref{thm-VMO-main}.
\end{proof}

\appendix
\section{}\label{appen}
  In this appendix, we will show $\mathcal{M}^{\beta_2}f$ is quasicontinuous on $Q_0$, whenever $f\in \bmob{\beta_1}(Q_0)$ for $0<\beta_1\leq \beta_2\leq n$ and  a finite open cube $Q_0$. Let us start with the following useful observations. The following Lemma is a little modification of \cite{ponce}*{Proposition~3.2} and also see \cite{BCRS}*{Remark~2.7}.
\begin{lemma}[\cite{ponce}]\label{a2}
 Let $Q$ be a finite cube of $\Rn$ and $1\leq p<\infty$, and $0<\beta\leq n$. If $f: Q\rightarrow \R$ is a quasicontinuous function and $\int_{Q}|f|^p\dc <\infty$, then there exists a sequence of continuous bounded functions $\{f_n\}_{n\geq 1}\subset C_b(Q)$ such that 
 \begin{align*}
\lim_{n\rightarrow\infty}\int_{Q}|f_n-f|^p\dc=0.
 \end{align*}
\end{lemma}

The next lemma says that if $f$ is continuous and bounded on a finite cube $Q$, then $\fdma f$ enjoys the same property on $Q$.

\begin{lemma}\label{a3}
Let $0<\beta\leq n$, and let $Q\subset \Rn$ be a finite cube. Suppose that $f$ is a real-valued function supported in $Q$ such that $f$ is continuous and bounded on $Q$. Then $\mathcal{M}^{\beta} f$ is continuous and bounded on $Q$.
\end{lemma}
\begin{proof}
It is easy to see that if $f$ is bounded, then $\mathcal{M}^{\beta} f$ is also bounded. Next, we can check that $\mathcal{M}^{\beta}f$ is lower semi-continuous. Hence, continuity of $\mathcal{M}^{\beta} f$ will follow by checking that $\mathcal{M}^{\beta}f$ is upper semi-continuous. This will be followed by simply using the continuity property of $f$. We also refer to \cite{mf-cont}*{Theorem~2.1} in this context.
\end{proof}

Now we are ready to state the result related to the quasicontinuity of $\mathcal{M}^{\beta_2}$.
\begin{theorem}\label{quasi}
Let $0 < \beta_1 \leq  \beta_2 \leq n$, and  $Q_0 \subset \mathbb{R}^n$ be a finite cube. Suppose $f \in \bmob{\beta_1}(Q_0)$ with $\supp(f) \subset Q_0$. Then $\mathcal{M}^{\beta_2} f$ is $\contb{\beta_1}$-quasicontinuous on $Q_0$.
\end{theorem}
\begin{proof}
Let us begin with $f \in \bmob{\beta_1}(Q_0)$. By applying the John–Nirenberg inequality, we deduce that $f\in L^p(\contb{\beta_1}, Q_0)$ for any $1<p<\infty$. By the definition of $L^p(\contb{\beta_1}, Q_0)$, it follows that $f$ is quasicontinuous. Moreover, by applying Lemma~\ref{a2}, there exists $\phi_j \in C_b(Q_0)$ such that
\begin{align}\label{p2}
  \lim_{j\rightarrow\infty}\int_{Q_0}|\phi_j-f|^p\dcb{\beta_1}=0.
\end{align}
Now, for $p\in\left(1, \infty\right)$, by applying \eqref{diff_max}, the fact that $\supp(f) \subset Q_0$, and Lemma~\ref{modf-lem}, we deduce that 
\begin{align*}   
\int_{Q_0}|\mathcal{M}^{\beta_2}\phi_j-\mathcal{M}^{\beta_2} f|^{p}\dcb{\beta_1} \leq \int_{Q_0}|\mathcal{M}^{\beta_2}\left(\phi_j-f\right)|^{p}\dcb{\beta_1} \leq C \left(\int_{Q_0}|\phi_j-f|^p\dcb{\beta_1}\right) .
\end{align*}
Hence, by applying \eqref{p2} above, we obtain that $\mathcal{M}^{\beta_2}\phi_j\rightarrow \mathcal{M}^{\beta_2} f$ in $L^{p}(\contb{\beta_1}, Q_0)$ as $j\rightarrow\infty$. Moreover, by Lemma~\ref{a3}, we have ${\mathcal{M}^{\beta_2}\phi_j}\subset L^{p}(\contb{\beta_1}, Q_0)$. Since this space is complete, it follows that the limit $\mathcal{M}^{\beta_2} f$ belongs to $L^{p}(\contb{\beta_1}, Q_0)$, and thus it is quasicontinuous on $Q_0$ by the definition of this space.
\end{proof}
    
\section*{Acknowledgments}
	R.~Basak is supported by the National Science and Technology Council of Taiwan under research grant numbers 113-2811-M-003-007/113-2811-M-003-039. Y.-W.~B.~Chen is supported by the National Science and Technology Council of Taiwan under research grant number 113-2811-M-002-027. P.~Roychowdhury is supported by the MUR-PRIN project no.~20227HX33Z, ``Pattern formation in nonlinear phenomena," granted by the European Union - Next Generation EU.
	
\end{document}